\documentclass{aptpubJuly17}

\authornames{Li-Xin Zhang} 
\shorttitle{Strong Limit theorems  for END random variables} 

\usepackage{mathrsfs}
\usepackage{bm}
\def\pr{\textsf{P}} 
\def\ep{\textsf{E}} 
\def\Sbep{\widehat{\mathbb E}} 
\def\cSbep{\widehat{\mathcal E}} 
\def\Capc{\mathbb V} 
\def\cCapc{\mathcal V} 
\def\outCapc{\mathbb V^{\ast}}
\def\outcCapc{\mathcal V^{\ast}}

\numberwithin{lemma}{section}
\numberwithin{theorem}{section}
\numberwithin{proposition}{section}
\numberwithin{assumption}{section}
\numberwithin{corollary}{section}
\numberwithin{definition}{section}
\numberwithin{remark}{section}
\numberwithin{example}{section}

\numberwithin{equation}{section}  

\begin{document}

\title{Strong limit theorems for extended independent  and extended negatively dependent random variables under non-linear expectations} 

\authorone[ Zhejiang University]{Li-Xin Zhang} 

\addressone{School  of Mathematical Sciences, Zhejiang University, Hangzhou 310027, P.R. China. Email:stazlx@zju.edu.cn} 

\begin{abstract}
Limit theorems for non-additive probabilities or non-linear expectations are challenging issues which have raised  progressive interest recently.  The purpose of this paper is to study the strong law of large numbers and the law of the iterated logarithm for a sequence of random variables in a sub-linear expectation space under a concept of  extended independence which is  much weaker and easier to verify  than the independence  proposed by Peng (2008b). We introduce a concept of extended  negative dependence which is an extension of this kind of weak  independence and the extended negative independence relative to classical probability appeared in recent literatures. Powerful tools as the moment inequality and Kolmogorov's exponential inequality are established for this kind of extended  negatively independent random variables, which improve those of Chen, Chen and Ng (2010) a lot.  And  the strong law of large numbers and the law of iterated logarithm are obtained by applying these inequalities.
\end{abstract}




\section{ Introduction.}\label{sect1}
\setcounter{equation}{0}

Non-additive probabilities and non-additive expectations  are useful tools for studying  uncertainties in statistics, measures of risk,
superhedging in finance and non-linear stochastic calculus, cf. Denis and  Martini (2006), Gilboa (1987),
Marinacci (1999),  Peng (1999, 2006, 2008a) etc.
This paper considers the  general sub-linear expectations and related non-additive probabilities generated by them. Under the frameworks of  the non-additive probability or non-linear expectation, the traditionary way for defining the independence is carried out  through the non-additive probability  by imitating the classical independence relative to the probability.  Under such frameworks, it is hard to study the limit theorems unless some additional  conditions (for example, the complete monotonicity of the non-additive probability) are assumed such that the non-additive probability is  somehow close to the additive one (c.f.  Maccheroni and Marinacci (2015), Ter\'an (2014)). To the best of my knowledge,
 Peng (2008b) is the first one to give a reasonable definition of the independence through the non-linear expectation.
Let $\{X_n;n\ge 1\}$ be a sequence of random variables in a sub-linear expectation space $(\Omega,\mathcal H, \Sbep)$. Peng's independence is that
 \begin{equation}\label{eq1}  \Sbep\left[ \psi(X_1,\cdots, X_n, X_{n+1})\right]=
\Sbep\left[ \Sbep\big[\psi(x_1,\cdots,x_n, X_{n+1})\big]\big|_{x_1=X_1,\cdots,x_n=X_n}\right],
\end{equation}
 for all $n\ge 2$ and any   $\psi\in C_{l,Lip}(\mathbb R_{n+1})$, where $ C_{l,Lip}(\mathbb R_{n+1})$ is the space of  local Lipschitz functions in $\mathbb R_{n+1}$. Under Peng's framework, many limit theorems have been being progressively established  very recently, including    the central limit theorem and weak law of large numbers (cf. Peng (2008b, 2010)),
 the law of the iterated algorithm (cf. Chen and Hu (2014), Zhang (2015a)), the small derivation and Chung's  law of the iterated logarithm (c.f. Zhang (2015b)),  the moment inequalities for the maximum partial sums (cf., Zhang (2016)).   Zhang (2016) gives the sufficient and necessary condition of the Kolomogov strong law of  large numbers.  For a sequence  of independent and identically distributed random variables $\{X_n;n\ge 1\}$, it is showed that the sufficient and necessary moment condition for the  strong law of large numbers to hold  is that the  Choquet integral    of $|X_1|$ is finite:
 \begin{equation}\label{eq2} C_{\Capc}(|X_1|)=\int_0^{\infty} \Capc(|X_1|\ge t)dt <\infty,
 \end{equation}
 where $\Capc$ is the upper capacity generated  by the sub-linear expectation $\Sbep$.

Recall that two random variables $X$ and $Y$ are independent relative to  a probability $P$ if and only if for any Borel functions $f$ and $g$,
$E_P[f(X)g(Y)]=E_P[f(X)]E_P[g(Y)]$
 whenever the expectations considered are finite. Another possible way to define the independence of $\{X_n;n\ge 1\}$ is that
 \begin{align}\label{eq3}
 &\Sbep\left[\psi_1(X_1,\ldots, X_k)\psi_2(X_{k+1},\ldots, X_n)\right]\nonumber\\
 = &
\Sbep\left[\psi_1(X_1,\ldots, X_k)\right]\Sbep\left[\psi_2(X_{k+1},\ldots, X_n)\right]
\end{align}
for all $n>k\ge 1$ and any $\psi_1\in C_{l,Lip}(\mathbb R_k)$ and $\psi_2\in C_{l,Lip}(\mathbb R_{n-k})$ such that the sub-linear expectations considered are finite.  If the independence is defined in this way, the functions $\psi_1$ and $\psi_2$ need to be limited in the class of non-negative functions, for otherwise we will conclude that $\Sbep[\cdot]=-\Sbep[-\cdot]$ and so $\Sbep$ will reduce to the linear expectation. It can showed that (\ref{eq1}) implies (\ref{eq3}). A more weaker independence is defined as the extended independence   in the sense that
\begin{equation}\label{eq4} \Sbep\left[\prod_{i=1}^n \psi_i(X_i)\right]=
\prod_{i=1}^n\Sbep\left[ \psi_i(X_i)\right], \;\forall\; n\ge 2, \forall\; 0\le \psi_i(x)\in C_{l,Lip}(\mathbb R).
\end{equation}
 This independence is much weaker than that of Peng and  easier to verify.  For the classical linear expectation, the above definitions of independence are equivalent. For the non-linear expectation, they are quite different. For example, Peng's independence has direction, i.e., that $Y$ is independent to $X$ does not imply that $X$ is independent to $Y$.  But the independence as in (\ref{eq4}) has no direction.
 One of the purposes of this paper is to show that, under this extended independence,    the sufficient and necessary moment condition for the Kolomogov strong law of  large numbers to hold is also that the  Choquet integral    of $|X_1|$ is finite. The proof of the sufficiency part is somewhat similar to that of Zhang (2016) after establishing good estimation of the tail capacity of partial sums of random variables.  Because we have not   ``the divergence part'' of the Borel-Cantelli Lemma and no information about the independence under the conjugate expectation $\cSbep$ or the conjugate capacity $\cCapc$, where $\cSbep[\cdot]=-\Sbep[-\cdot]$ and $\cCapc(A)=1-\Capc(A^c)$,   proving the necessary part is a challenging work.

By replacing the function space $C_{l,lip}(\mathbb R)$ with the family of all Borel measurable functions, Chen, Wu and Li (2013) considered  random variables which are independent  in sense of (\ref{eq4}) under a upper expectation $\Sbep[\cdot]$ being defined by
$$  \Sbep[X]=\sup_{P\in \mathscr{P}} E_P[X],$$
where $\mathscr{P}$ is a family of probability measures defined on a measurable space $(\Omega, \mathcal{F})$. The strong law of large numbers was proved under finite $(1+\alpha)$-th moments $(\alpha>0)$ which is much stringer than (\ref{eq2}) when the random variables are identically distributed.  Note that, if $\{X_n;n\ge 1\}$ are independent relative   to each $P\in \mathscr{P}$, then we will have
\begin{equation} \label{eq5} \Sbep\left[\prod_{i=1}^n \psi_i(X_i)\right]\le
\prod_{i=1}^n\Sbep\left[ \psi_i(X_i)\right], \;\forall\; n\ge 2, \forall\; 0\le \psi_i(x)\in C_{l,Lip}(\mathbb R).
\end{equation}
But in general, the equality will not hold. So the random variables  may not be independent under $\Sbep$.
A simple example is that $$(X_1,X_2)\sim P_{\sigma}\in \mathscr{P}=\{N(0,\sigma^2)\otimes N(0,\sigma^{-2}):
1/2\le \sigma\le 2\}$$
for which $\sup\limits_{\sigma}E_{\sigma}\left[X_1^2X_2^2\right]=1$, $\sup\limits_{\sigma}E_{\sigma}\left[X_1^2\right]\sup\limits_{\sigma}E_{\sigma}\left[X_1^2\right]=16.$
It is of important interest to study the limit theorems for random variables satisfying the property (\ref{eq5}).

The property (\ref{eq5})  is very close to that of negatively dependent random variables.  The concept of negative dependence relative to the classical probability has been extensively investigated since it appeared in Lehmann (1966). Various generalization of the concept of negative dependence and related limit theorems have been studied in literatures. One can refer to  Joag-Dev  and Proschan (1983), Newman (1984), Matula (1992), Su et al  (1997), Shao  and Su  (1999), Shao (2000), Zhang (2001a, 2001b) etc.  As a new extension, the concept of  extended negative  dependence was proposed in Liu (2009) and further promoted in Chen, Chen and Ng (2010). A sequence of random variables is said to be extended negatively dependent  if the tails of its finite-dimensional distributions in the lower-left and upper-right corners are dominated by a multiple of the tails of the corresponding finite-dimensional distributions of a sequence of independent random variables with the same marginal distributions. The strong law of large numbers was   established by Chen, Chen and Ng (2010). However, for the extended negatively  dependent random variables, besides the type of the law of large numbers, very little is known about other kinds of fine limit theorems such as the central limit theorem and the law of the iterated logarithm.   In this paper, we will introduce a concept of extended  negative dependence under the sub-linear expectation which is weaker than the extended independence as defined in (\ref{eq4}) and is an extension of the extended negative  dependence relative to the classical probability. The strong law of large numbers will also be established for    extended  negatively dependent random variables. The result of Chen, Chen and Ng (2010) is extended and improved.

To establish the strong law of large numbers, some key inequalities for  the tails of the capacities of the sums of extended negatively  dependent random variables in the general sub-linear expectation spaces are obtained, including the moment inequalities and the Kolmogorov type exponential inequalities. These inequalities also improve those established by Chen, Chen and Ng (2010) for extended negatively dependent random variables relative to a classical probability, as well as those for independent random variables in a sub-linear expectation space.  They  may be useful tools for studying other limit theorems. We also establish  the law of the iterated logarithm by applying the exponential inequalities. And as a corollary, the law of the iterated logarithm for extended negatively dependent random variables on a probability space is obtained.
    In the next section, we give some notations under the sub-linear expectations.
    In Section~\ref{sectIneq}, we will establish the exponential inequalities.   The law of large numbers is given in Section~\ref{sectLLN}. In the last section we consider the law of the iterated logarithm.

\section{Basic Settings}\label{sectBasic}
\setcounter{equation}{0}

We use the framework and notations of Peng (2008b). Let  $(\Omega,\mathcal F)$
 be a given measurable space  and let $\mathscr{H}$ be a linear space of real functions
defined on $(\Omega,\mathcal F)$ such that if $X_1,\ldots, X_n \in \mathscr{H}$  then $\varphi(X_1,\ldots,X_n)\in \mathscr{H}$ for each
$\varphi\in C_{l,Lip}(\mathbb R_n)$,  where $C_{l,Lip}(\mathbb R_n)$ denotes the linear space of (local Lipschitz)
functions $\varphi$ satisfying
\begin{eqnarray*} & |\varphi(\bm x) - \varphi(\bm y)| \le  C(1 + |\bm x|^m + |\bm y|^m)|\bm x- \bm y|, \;\; \forall \bm x, \bm y \in \mathbb R_n,&\\
& \text {for some }  C > 0, m \in \mathbb  N \text{ depending on } \varphi. &
\end{eqnarray*}
$\mathscr{H}$ is considered as a space of ``random variables''. In this case we denote $X\in \mathscr{H}$.
We also denote $C_{b,Lip}(\mathbb R_n)$ to be the bounded Lipschitz functions $\psi(x)$ satisfying
  \begin{eqnarray*} & |\varphi(\bm x)|\le C,\;\; |\varphi(\bm x) - \varphi(\bm y)| \le  C|\bm x- \bm y|, \;\; \forall \bm x, \bm y \in \mathbb R_n,&\\
& \text {for some }  C > 0,  \text{ depending on } \varphi. &
\end{eqnarray*}
\begin{definition}\label{def1.1} A  sub-linear expectation $\Sbep$ on $\mathscr{H}$  is a function $\Sbep: \mathscr{H}\to \overline{\mathbb R}$ satisfying the following properties: for all $X, Y \in \mathscr H$, we have
\begin{description}
  \item[\rm (a)]  Monotonicity: If $X \ge  Y$ then $\Sbep [X]\ge \Sbep [Y]$;
\item[\rm (b)] Constant preserving: $\Sbep [c] = c$;
\item[\rm (c)] Sub-additivity: $\Sbep[X+Y]\le \Sbep [X] +\Sbep [Y ]$ whenever $\Sbep [X] +\Sbep [Y ]$ is not of the form $+\infty-\infty$ or $-\infty+\infty$;
\item[\rm (d)] Positive homogeneity: $\Sbep [\lambda X] = \lambda \Sbep  [X]$, $\lambda\ge 0$.
 \end{description}
 Here $\overline{\mathbb R}=[-\infty, \infty]$. The triple $(\Omega, \mathscr{H}, \Sbep)$ is called a sub-linear expectation space. Give a sub-linear expectation $\Sbep $, let us denote the conjugate expectation $\cSbep$of $\Sbep$ by
$$ \cSbep[X]:=-\Sbep[-X], \;\; \forall X\in \mathscr{H}. $$
\end{definition}

From the definition, it is easily shown that    $\cSbep[X]\le \Sbep[X]$, $\Sbep[X+c]= \Sbep[X]+c$ and $\Sbep[X-Y]\ge \Sbep[X]-\Sbep[Y]$ for all
$X, Y\in \mathscr{H}$ with $\Sbep[Y]$ being finite. Further, if $\Sbep[|X|]$ is finite, then $\cSbep[X]$ and $\Sbep[X]$ are both finite.

Next, we consider the capacities corresponding to the sub-linear expectations.
Let $\mathcal G\subset\mathcal F$. A function $V:\mathcal G\to [0,1]$ is called a capacity if
$$ V(\emptyset)=0, \;V(\Omega)=1 \; \text{ and } V(A)\le V(B)\;\; \forall\; A\subset B, \; A,B\in \mathcal G. $$
It is called to be sub-additive if $V(A\bigcup B)\le V(A)+V(B)$ for all $A,B\in \mathcal G$  with $A\bigcup B\in \mathcal G$.

In the sub-linear space $(\Omega, \mathscr{H}, \Sbep)$, we denote a pair $(\Capc,\cCapc)$ of capacities by
$$ \Capc(A):=\inf\{\Sbep[\xi]: I_A\le \xi, \xi\in\mathscr{H}\}, \;\; \cCapc(A):= 1-\Capc(A^c),\;\; \forall A\in \mathcal F, $$
where $A^c$  is the complement set of $A$.
Then
\begin{equation}\label{eq1.3} \begin{matrix}
&\Capc(A):=\Sbep[I_A], \;\; \cCapc(A):= \cSbep[I_A],\;\; \text{ if } I_A\in \mathscr H\\
&\Sbep[f]\le \Capc(A)\le \Sbep[g], \;\;\cSbep[f]\le \cCapc(A) \le \cSbep[g],\;\;
\text{ if } f\le I_A\le g, f,g \in \mathscr{H}.
\end{matrix}
\end{equation}
It is obvious that $\Capc$ is sub-additive. But $\cCapc$ and $\cSbep$ are not. However, we have
\begin{equation}\label{eq1.4}
  \cCapc(A\bigcup B)\le \cCapc(A)+\Capc(B) \;\;\text{ and }\;\; \cSbep[X+Y]\le \cSbep[X]+\Sbep[Y]
\end{equation}
due to the fact that $\Capc(A^c\bigcap B^c)\ge \Capc(A^c)-\Capc(B)$ and $\Sbep[-X-Y]\ge \Sbep[-X]-\Sbep[Y]$.

Also, we define the  Choquet integrals/expecations $(C_{\Capc},C_{\cCapc})$  by
$$ C_V[X]=\int_0^{\infty} V(X\ge t)dt +\int_{-\infty}^0\left[V(X\ge t)-1\right]dt $$
with $V$ being replaced by $\Capc$ and $\cCapc$ respectively. It is obvious that
$$ C_{\Capc}(|X|)\le 1+ \Sbep[|X|^{1+\alpha}]\int_1^{\infty} t^{-1-\alpha}dt \le 1+\alpha^{-1} \Sbep[|X|^{1+\alpha}]. $$
 Also,
it can be verified that, if $\lim_{c\to \infty}\Sbep[(|X|-c)^+]=0$, then
$ \Sbep[|X|]\le C_{\Capc}(|X|)$ (c.f. Lemma 3.9 of Zhang (2016)).

\bigskip

The concept of independence and identical distribution is introduced by Peng (2006,2008b).

\begin{definition} ({\em Peng (2006, 2008b)})

\begin{description}
  \item[ \rm (i)] ({\em Identical distribution}) Let $\bm X_1$ and $\bm X_2$ be two $n$-dimensional random vectors defined
respectively in sub-linear expectation spaces $(\Omega_1, \mathscr{H}_1, \Sbep_1)$
  and $(\Omega_2, \mathscr{H}_2, \Sbep_2)$. They are called identically distributed, denoted by $\bm X_1\overset{d}= \bm X_2$  if
$$ \Sbep_1[\varphi(\bm X_1)]=\Sbep_2[\varphi(\bm X_2)], \;\; \forall \varphi\in C_{l,Lip}(\mathbb R_n), $$
whenever the sub-expectations are finite. A sequence $\{X_n;n\ge 1\}$ of random variables is said to be identically distributed if $X_i\overset{d}= X_1$ for each $i\ge 1$.
\item[\rm (ii)] ({\em Independence})   In a sub-linear expectation space  $(\Omega, \mathscr{H}, \Sbep)$, a random vector $\bm Y =
(Y_1, \ldots, Y_n)$, $Y_i \in \mathscr{H}$ is said to be independent to another random vector $\bm X =
(X_1, \ldots, X_m)$ , $X_i \in \mathscr{H}$ under $\Sbep$  if for each test function $\varphi\in C_{l,Lip}(\mathbb R_m \times \mathbb R_n)$
we have
$ \Sbep [\varphi(\bm X, \bm Y )] = \Sbep \big[\Sbep[\varphi(\bm x, \bm Y )]\big|_{\bm x=\bm X}\big],$
whenever $\overline{\varphi}(\bm x):=\Sbep\left[|\varphi(\bm x, \bm Y )|\right]<\infty$ for all $\bm x$ and
 $\Sbep\left[|\overline{\varphi}(\bm X)|\right]<\infty$.
 \item[\rm (iii)] ({\em Independent random variables}) A sequence of random variables $\{X_n; n\ge 1\}$
 is said to be independent, if
 $X_{i+1}$ is independent to $(X_{1},\ldots, X_i)$ for each $i\ge 1$.
 \end{description}
\end{definition}

\begin{definition}  \label{WeakInd}
  ({\em Extended  Independence}) A sequence of random variables $\{X_n; n\ge 1\}$
 is said to be extended independent, if
\begin{equation}\label{eqWeakInd} \mathbb{E}\left[\prod_{i=1}^n \psi_i(X_i)\right]=
\prod_{i=1}^n\mathbb{E}\left[ \psi_i(X_i)\right], \;\forall\; n\ge 2, \forall\; 0\le \psi_i(x)\in C_{l,Lip}(\mathbb R).
\end{equation}
\end{definition}

It can be showed that the independence implies the extended independence. It shall be noted that the extended independence of $\{X_n; n\ge 1\}$ under $\Sbep$ does not imply the extended independence under $\cSbep$.  The independence in sense of (\ref{eqWeakInd}) was proposed in Chen, Wu and Li (2013). But their function space of $\psi_i$s is assumed to be the family of all non-negative Borel functions. Here we use the function space the same as Peng's. The function space can also be limited to $C_{b,Lip}(\mathbb R)$.

  Recall that a sequence of random variables $\{Y_n;n\ge 1\}$  on a probability space $(\Omega, \mathcal F, \pr)$ are called to be lower  extended negatively dependent (LEND)  if there is some dominating constant $K \ge  1$ such that, for all $x_i$, $i = 1, 2,\ldots$,
\begin{equation}\label{eqENDlower} \pr\left(\bigcap_{i=1}^n \left\{Y_i\le x_i\right\}\right)\le K \prod_{i=1}^n \pr\left(Y_i\le x_i\right), \;\; \forall\; n,
\end{equation}
and they are called upper  extended negatively dependent (UEND)  if   for all $x_i$, $i = 1, 2,\ldots$,
\begin{equation}\label{eqENDupper}\pr\left(\bigcap_{i=k=1}^n \left\{Y_i> x_i\right\}\right)\le K \prod_{i=1}^n \pr\left(Y_i> x_i\right), \;\; \forall\; n.
\end{equation}
They are called extended negatively dependent (END) if they are both LEND and UEND (cf., Liu (2009)). In the case $K=1$ the notion of END random variables reduces to the well known
notion of so-called negatively dependent (ND) random variables which was
introduced by Lehmann (1966) (cf. also  Block et al. (1982), Joag-Dev and Proschan (1983) etc).    It is showed that if $\{Y_n;n\ge 1\}$ are  upper (resp. lower)  extended negatively dependent, and the functions $g_i\ge 0$, $i = 1,2,\ldots$, are all non-decreasing (resp. all non-increasing), then
\begin{equation}\label{eqEND}
\ep\left[\prod_{i=1}^ng_i(Y_i)\right]\le K\prod_{i=1}^n\ep\left[g_i(Y_i)\right], \; n\ge 1,
\end{equation}
(cf., Chen, Chen and Ng (2010)).
Motivated by the above property (\ref{eqEND}) and Definition \ref{WeakInd},    we introduce  a concept of extended negative dependence under the sub-linear expectation.

 \begin{definition}\label{Def2.4}
   ({\em Extended negative dependence})  In a sub-linear expectation space  $(\Omega, \mathscr{H}, \Sbep)$,  random variables $\{X_n;n\ge 1\}$   are called to be upper (resp. lower) extended negatively dependent  if there is some   dominating constant $K \ge 1$ such that \begin{equation}\label{eqWND}
\Sbep\left[\prod_{i=1}^ng_i(X_i)\right]\le K\prod_{i=1}^n\Sbep\left[g_i(X_i)\right], n\ge 1,
\end{equation}
 whenever the non-negative functions $g_i\in C_{b,Lip}(\mathbb R)$, $i = 1,2,\ldots$, are all non-decreasing (resp. all non-increasing). They are called extended negatively dependent if they both upper extended negatively dependent and lower extended negatively dependent.
\end{definition}
It is obvious that, if $\{X_n;n\ge 1\}$ is a sequence of  extended independent random variables and $f_1(x),f_2(x),\ldots\in C_{l,Lip}(\mathbb R)$,
 then $\{f_n(X_n);n\ge 1\}$ is  also a sequence of extended independent random variables, and they are   extended negatively dependent with $K=1$;  if $\{X_n;n\ge 1\}$ is  a sequence of upper extended negatively dependent random variables and $f_1(x),f_2(x),\ldots\in C_{l,Lip}(\mathbb R)$
 are all non-decreasing (resp. all non-increasing) functions, then $\{f_n(X_n);n\ge 1\}$ is  also a sequence of upper (resp. lower) extended negatively dependent random variables.

\begin{example}
Let $(\Omega,\mathcal F)$ be a measurable space, $\mathscr{P}$ be a family  of probability measures on it and $\{X_n;n\ge 1\}$ be a sequence of random variables. Define a upper expectation $\Sbep[\cdot]$ by
$$ \Sbep[X]=\sup_{P\in \mathscr{P}}E_P[X]. $$
Then $\Sbep$ is a sub-linear expectation. If $\{X_n;n\ge 1\}$ are extended negatively dependent in the sense  of (\ref{eqENDupper}) and (\ref{eqENDlower}) relative to each $P\in \mathscr{P}$ with the same  dominating constant $K$, then  they are extended negatively dependent under $\Sbep$.
\end{example}
We will establish the exponential inequalities, the law of large numbers and the law of the iterated logarithm for this kind of extended independent random variables.

\section{ Exponential inequalities}\label{sectIneq}
\setcounter{equation}{0}
In this section, we are going to establish some key inequalities   for the sums of extended negatively  dependent random variables, including  moment inequalities and   the  exponential inequalities. These inequalities improve Lemmas 2.5 and 2.6 of Chen, Chen and Ng (2010).
Let $\{X_1,\ldots, X_n\}$ be a sequence  of   random variables
in $(\Omega, \mathscr{H}, \Sbep)$. Set $S_n=\sum_{k=1}^n X_k$,  $B_n=\sum_{k=1}^n \Sbep[X_k^2]$ and $M_{n,p}=\sum_{k=1}^n\Sbep[|X_k|^p]$, $M_{n,p,+}=\sum_{k=1}^n\Sbep[(X_k^+)^p]$, $p\ge 2$.

\begin{theorem}\label{thIneq2}    Let $\{X_1,\ldots, X_n\}$ be a sequence  of upper extended  negatively dependent random variables
in $(\Omega, \mathscr{H}, \Sbep)$ with $\Sbep[X_k]\le 0$. Then
\begin{description}
  \item[\rm (a)]
  For all $x,y>0$,
\begin{align}\label{eqthIneq2.1}
\Capc\left(S_n\ge x\right)\le &\Capc\left(\max_{k\le n}X_k\ge y\right) \nonumber \\
&+ K\exp\left\{-\frac{x^2}{2(xy+B_n)}\Big(1+\frac{2}{3}\ln \big(1+\frac{xy}{B_n}\big)\Big)\right\};
\end{align}
\item[\rm (b)]  For any $p\ge 2$, there exists a constant $C_p\ge 1$ such that for all $x>0$ and $0<\delta\le 1$,
\begin{equation}\label{eqthIneq2.2}
\Capc\left(S_n\ge x\right)\le C_p \delta^{-2p} K\frac{M_{n,p,+}}{x^p}
+K\exp\left\{-\frac{x^2}{2 B_n (1+\delta)}\right\},
\end{equation}
\item[\rm (c)] We have for $x>0$, $r>0$ and $p\ge 2$,
 \begin{equation}\label{eqthIneq2.3ad}
   \Capc\left(S_n^+\ge x\right)\le   \Capc\big(\max_{k\le n}X_k^+\ge \frac{x}{r}\big)
+K e^r\left( \frac{r B_n}{r B_n+x^2}\right)^{r},
\end{equation}
\begin{align}\label{eqthIneq2.3}C_{\Capc}\left[(S_n^+)^p\right]
\le  & p^p C_{\Capc}\Big[\big(\max_{k\le n}X_k^+\big)^p\Big]+ C_p B_n^{p/2} \nonumber \\
\le &
p^p \sum_{k=1}^n C_{\Capc}\Big[(X_k^+)^p\Big]+ C_p B_n^{p/2}.
\end{align}
In particular,
\begin{equation}\label{eqthIneq2.2ad}
\Capc\left(S_n\ge x\right)\le (1+Ke)  \frac{B_n}{x^2}, \;\; \forall x>0.
\end{equation}
\end{description}
 \end{theorem}

\proof  Let $Y_k=X_k\wedge y$, $T_n=\sum_{k=1}^n Y_k$.  Then $X_k-Y_k=(X_k-y)^+\ge 0$ and
 $\Sbep[Y_k]\le \Sbep[X_k]\le 0$.
Note that $\varphi(x)=:e^{t(x\wedge y)}$ is a bounded non-decreasing function and belongs to $C_{b,Lip}(\mathbb R)$ since $0\le \varphi^{\prime}(x)\le t e^{ty}$ if $t>0$.
  It follows that for any $t>0$,
\begin{align*}
 \Capc\left(S_n\ge x\right)\le &  \Capc\big(\max_{k\le n}X_k\ge y\big)+
 \Capc\left(T_n\ge x\right),\\
\Capc\left(T_n\ge x\right)  \le &
 e^{-tx}\Sbep[e^{t T_n}]
 \le
 e^{-tx}K\prod_{k=1}^n \Sbep[e^{t Y_k}],
 \end{align*}
be the definition of the upper extended negative dependence.
The remainder of the proof is the similar to that of Zhang (2015a). For the completeness of this paper, we also present it here.

Note
$$e^{tY_k}=1+ tY_k+\frac{e^{tY_k}-1-t Y_k}{Y_k^2}Y_k^2\le 1 +tY_k+\frac{e^{ty}-1-t y}{y^2}Y_k^2.  $$
We have
$$ \Sbep[e^{t Y_k}]\le 1+\frac{e^{ty}-1-t y}{y^2} \Sbep[Y_k^2]\le\exp\left\{\frac{e^{ty}-1-t y}{y^2}\Sbep[X_k^2]\right\}. $$
Choosing $t=\frac{1}{y}\ln \big(1+\frac{xy}{B_n}\big)$ yields
\begin{align}\label{eqproofth2.1}
\Capc\left(T_n\ge x\right)  \le &K e^{-tx}\exp\left\{\frac{e^{ty}-1-t y}{y^2}B_n\right\}\nonumber \\
=&K\exp\left\{\frac{x}{y}-\frac{x}{y}\Big(\frac{B_n}{xy}+1\Big)\ln\Big(1+\frac{xy}{B_n}\Big)\right\}.
 \end{align}
  Applying the elementary inequality
 $$ \ln (1+t)\ge \frac{t}{1+t}+\frac{t^2}{2(1+t)^2}\big(1+\frac{2}{3} \ln (1+t)\big)$$
 yields
 $$ \Big(\frac{B_n}{xy}+1\Big)\ln\Big(1+\frac{xy}{B_n}\Big)
 \ge 1+\frac{xy}{2(xy+B_n)}\Big(1+\frac{2}{3}\ln\big(1+\frac{xy}{B_n}\big)\Big). $$
 (\ref{eqthIneq2.1}) is proved.

 Next we show (b). If $xy\le \delta B_n$, then
 $$\frac{x^2}{2(xy+B_n)}\Big(1+\frac{2}{3}\ln\big(1+\frac{xy}{B_n}\big)\Big)\ge \frac{x^2}{2 B_n(1+\delta)}. $$
If  $xy\ge \delta B_n$,
then
 $$\frac{x^2}{2(xy+B_n)}\Big(1+\frac{2}{3}\ln\big(1+\frac{xy}{B_n}\big)\Big)\ge \frac{x}{2(1+1/\delta)y}. $$
 It follows that
\begin{equation}\label{eqproofth2.2} \Capc\left(T_n\ge x\right)\le K\exp\left\{-\frac{x^2}{2 B_n (1+\delta)}\right\}
 +K \exp\left\{-\frac{x}{2 (1+1/\delta)y}\right\}
 \end{equation}
 by (\ref{eqproofth2.1}). Let
 $$\beta(x)=\beta_p(x)=\frac{1}{x^p}\sum_{k=1}^n\Sbep[(X_k^+)^p], $$
 and choose
 $$ \rho=1\wedge \frac{1}{2(1+1/\delta)\delta \log (1/\beta(x))}, \;\; y=\rho\delta x. $$
 Then by (\ref{eqproofth2.2}),
 \begin{align*}
  &\Capc\big(S_n\ge (1+2\delta)x\big)\le  \Capc\big(T_n\ge  x\big)+\Capc\big(\sum_{i=1}^n(X_i-\rho\delta x)^+\ge 2\delta x\big) \\
  \le &K\exp\left\{-\frac{x^2}{2 B_n (1+\delta)}\right\}+K\beta(x)+\Capc\big(\max_{i\le n} X_i\ge \delta x\big)+\Capc\big(\sum_{i=1}^n(X_i-\rho\delta x)^+\wedge (\delta x) \ge 2\delta x\big).
  \end{align*}
  It is obvious that
  $$ \Capc\big(\max_{i\le n} X_i\ge \delta x\big)\le \delta^{-p} \beta(x). $$
  On the other hand,
  \begin{align*}
  &\Capc\big(\sum_{i=1}^n(X_i-\rho\delta x)^+\wedge (\delta x) \ge 2\delta x\big)=
 \Capc\left(\sum_{i=1}^n\left[\Big(\frac{X_i}{\delta x}-\rho\big)^+\wedge 1\right] \ge 2\right)\\
 \le & e^{-2t}\Sbep\exp\left\{t\sum_{i=1}^n\left[\Big(\frac{X_i}{\delta x}-\rho\big)^+\wedge 1\right]\right\}
 \le e^{-2t}K\prod_{i=1}^n \Sbep\exp\left\{t\left[\Big(\frac{X_i}{\delta x}-\rho\big)^+\wedge 1\right]\right\}\\
 \le & e^{-2t}K\prod_{i=1}^n\left[ 1+e^t \Capc(X_i\ge \rho\delta x)\right]
 \le K\exp\left\{ -2t +e^t \sum_{i=1}^n \Capc(X_i\ge \rho\delta x)\right\},
 \end{align*}
 where the second inequality is due to the upper extended negative dependence.   Assume  $\beta(x)<1$. Let $e^t \sum_{i=1}^n \Capc(X_i\ge \rho\delta x)=2$ (while, if $\sum_{i=1}^n \Capc(X_i\ge \rho\delta x)=0$, we let $t\to \infty$). We obtain
\begin{align*}
&  \Capc\big(\sum_{i=1}^n(X_i-\rho\delta x)^+\wedge (\delta x) \ge 2\delta x\big)
\le   Ke^2 \left(\frac{1}{2}\sum_{i=1}^n \Capc(X_i\ge \rho\delta x)\right)^2\\
\le & Ke^2 \left(\frac{\beta(x)}{2(\delta \rho)^p}\right)^2
=K e^22^{-2}\delta^{-2p} (2(\delta+1))^{2p}\beta^2(x)\left(\log\frac{1}{\beta(x)}\right)^{2p}
\le    C_p\delta^{-2p}\beta(x),
\end{align*}
where the last inequality is due to the fact that $(\log 1/t)^{2p}\le C_p /t$ ($0<t<1$). It follows that
\begin{align*}
  \Capc\big(S_n\ge (1+2\delta)x\big)\le K\exp\left\{-\frac{x^2}{2 B_n (1+\delta)}\right\}+C_p K \delta^{-2p}\beta(x).
  \end{align*}
 If $\beta(x)\ge 1$, then the above inequality is obvious. Now letting $z= (1+2\delta)x $ and $\delta^{\prime}=(1+\delta)(1+2\delta)^2-1$ yields
 \begin{align*}
  \Capc\big(S_n\ge z\big)\le K\exp\left\{-\frac{z^2}{2 B_n (1+\delta^{\prime})}\right\}+C_p K (\delta^{\prime})^{-2p}\beta(z).
  \end{align*}
 (b) is proved.

Finally, we consider (c).
Putting $y=x/r$ in (\ref{eqproofth2.1}), we obtain (\ref{eqthIneq2.3ad}).
Note that
$$ C_{\Capc}\big[(X^+)^p\big]= \int_0^{\infty}\Capc\big(X^p>x)dx= \int_0^{\infty}px^{p-1}\Capc\big(X>x)dx. $$
Then putting   $r=p>p/2$, multiplying  both sides of  (\ref{eqthIneq2.3ad})  by $px^{p-1}$ and
  integrating on the positive half-line, we conclude (\ref{eqthIneq2.3}).  \endproof

  \section{The law of Large numbers}\label{sectLLN}
\setcounter{equation}{0}
For a sequence $\{X_n;n\ge 1\}$ of  random variables in the sub-linear expectation space $(\Omega, \mathscr{H}, \Sbep)$,
we denote $S_n=\sum_{k=1}^n X_k$, $S_0=0$. We first consider the weak law of large numbers.
 \begin{theorem} \label{thWLLN1}
  \begin{description}
  \item[\rm (a)]   Suppose that $X_1, X_2, \ldots$ are identically distributed and  extended negatively  dependent with    $\lim\limits_{c\to \infty} \Sbep\left[(|X_1|-c)^+\right]=0$.
   Then for any $\epsilon>0$,
\begin{equation}\label{eqthWLLN1.1}
\cCapc\left( \cSbep[X_1]-\epsilon\le  \frac{S_n}{n}\le  \Sbep[X_1]+\epsilon\right)\to 1.
\end{equation}
\item[\rm (b)]   Suppose that $X_1, X_2, \ldots$ are  identically distributed and extended independent with    $\lim\limits_{c\to \infty} \Sbep\left[(|X_1|-c)^+\right]=0$. Then for any $\epsilon>0$,
  \begin{equation}\label{eqthWLLN1.2}
\Capc\left(\Big|\frac{S_n}{n}-\Sbep[X_1]\Big|\le \epsilon\right)\to 1
\end{equation}
and
\begin{equation}\label{eqthWLLN1.3}
\Capc\left(\Big|\frac{S_n}{n}-\cSbep[X_1]\Big|\le \epsilon\right)\to 1.
\end{equation}
\end{description}
\end{theorem}
We conjuncture that for any point $p\in  \big[\cSbep[X_1],\Sbep[X_1]\big]$,
 $\Capc\left(\Big|\frac{S_n}{n}-p\Big|\le \epsilon\right)\to 1.$

\proof Define
\begin{equation}\label{eqproofthWLLN1.1} f_c(x)= (-c)\vee (x\wedge c), \;\; \widehat{f}_c(x)=x-f_c(x)
\end{equation}
 and $\overline{X}_j=f_c(X_j)$,  $j=1,2,\ldots$.
 Then $f_c(\cdot)\in C_{l,Lip}(\mathbb R)$, and $\overline{X}_j$, $j=1,2,\ldots$ are extended negatively dependent (resp. extended independent) identically distributed random variables. It is easily seen that $\Sbep[f_c(X_1)]\to \Sbep[ X_1]$, $\cSbep[f_c(X_1)]\to \cSbep[ X_1]$ as $c\to +\infty$, and
\begin{align*} \sup_n\Capc\left(\Big|S_n-\sum_{j=1}^n \overline{X}_j\Big|\ge \epsilon n\right)
\le & \sup_n \frac{1}{\epsilon n}  \sum_{j=1}^n \Sbep|f_c(X_j)|  \\
\le & \frac{1}{\epsilon} \Sbep(|X_1|-c)^+\to 0 \text{ as } c\to \infty.
\end{align*}
So, it is sufficient to consider $\{\overline{X}_n; n\ge 1\}$ and then without loss of generality we can assume that $X_n$ is bounded by a constant $c>0$. By (\ref{eqthIneq2.2ad}),
 \begin{equation} \label{eqthproofWLLN.1} \Capc\left(\frac{S_n}{n}-\Sbep[X_1]\ge \epsilon\right)\le C \frac{\sum_{j=1}^n \Sbep\left[\big(X_j-\Sbep[X_j]\big)^2\right]}{n^2\epsilon^2}\le C \frac{n c^2  }{n^2\epsilon^2}\to 0,
 \end{equation}
and similarly,
$$ \Capc\left(\frac{S_n}{n}-\cSbep[X_1]\le -\epsilon\right)=  \Capc\left(\frac{-S_n}{n}-\Sbep[-X_1]\ge \epsilon\right)\to 0. $$
(\ref{eqthWLLN1.1}) is proved.

Now, suppose that $X_1, X_2, \ldots,$ are   extended independent. By noting (\ref{eqthWLLN1.1}), for (\ref{eqthWLLN1.2}) and (\ref{eqthWLLN1.3}) it is sufficient to show that
\begin{equation}\label{eqthproofWLLN.2}
\Capc\left(\frac{S_n}{n}-\Sbep[X_1]\ge -\epsilon\right)\to 1.
\end{equation}
It is easily  seen that   (\ref{eqthproofWLLN.2}) is equivalent to
 \begin{equation}\label{eqthproofWLLN.2ad}
\cCapc\left(\frac{\sum_{k=1}^n (-X_j-\cSbep[-X_j])}{n}\ge  \epsilon\right)\to 0.
\end{equation}
However, we have no inequality to estimate the lower capacity $\cCapc(\cdot)$ such that (\ref{eqthproofWLLN.2ad}) can be verified.
We shall now introduce a more technique method to show  (\ref{eqthproofWLLN.2}).

For any $0<\delta<\epsilon$ and $t>0$, we have
\begin{align*}
&I\left\{\frac{S_n}{n}-\Sbep[X_1]\ge -\epsilon\right\}\\
\ge &   e^{-t\delta}\left(\exp\left\{t\frac{S_n-n\Sbep[X_1]}{n}\right\}-e^{-t\epsilon}\right)I\left\{\frac{S_n-n\Sbep[X_1]}{n}\le  \delta\right\}\\
= & e^{-t\delta}\left(\exp\left\{t\frac{S_n-n\Sbep[X_1]}{n}\right\}-e^{-t\epsilon}\right)\\
&-e^{-t\delta}\left(\exp\left\{t\frac{S_n-n\Sbep[X_1]}{n}\right\}-
e^{-t\epsilon}\right)I\left\{\frac{S_n-n\Sbep[X_1]}{n}>  \delta\right\}\\
\ge & e^{-t\delta}\left(\prod_{j=1}^n \exp\left\{t\frac{X_j-\Sbep[X_j]}{n}\right\}-e^{-t\epsilon}\right)
-e^{-t\delta}e^{2tc}I\left\{\frac{S_n-n\Sbep[X_1]}{n}>  \delta\right\}.
\end{align*}
It follows that
\begin{align*}
&\Capc\left(\frac{S_n}{n}-\Sbep[X_1]\ge -\epsilon\right)\\
\ge & e^{-t\delta}\left(\Sbep\left[\prod_{j=1}^n \exp\left\{t\frac{X_j-\Sbep[X_j]}{n}\right\}\right]-e^{-t\epsilon}\right)
-e^{-t\delta}e^{2tc}\Capc\left(\frac{S_n-n\Sbep[X_1]}{n}>  \delta\right).
\end{align*}
By (\ref{eqthproofWLLN.1}), the second term will goes to zero as $n\to\infty$. By the extended independence and the fact that $e^x\ge 1+x$,
\begin{align*}
\Sbep\left[\prod_{j=1}^n \exp\left\{t\frac{X_j-\Sbep[X_j]}{n}\right\}\right]
=&\prod_{j=1}^n \Sbep\left[\exp\left\{t\frac{X_j-\Sbep[X_j]}{n}\right\}\right]\\
\ge &\prod_{j=1}^n \Sbep\left[ t\frac{X_j-\Sbep[X_j]}{n}+1\right]=1.
\end{align*}
It follows that
$$ \liminf_{n\to \infty} \Capc\left(\frac{S_n}{n}-\Sbep[X_1]\ge -\epsilon\right)
\ge e^{-t\delta}\left(1-e^{-t\epsilon}\right). $$
Letting $\delta\to 0$  and then $t\to \infty$ yields (\ref{eqthproofWLLN.2}). The proof is completed.  \endproof

\bigskip
Before we give the strong laws of large numbers,  we need some more notations about the sub-linear expectations and capacities.

\begin{definition}\label{def3.1}
\begin{description}
\item{\rm (I)} A sub-linear expectation $\Sbep: \mathscr{H}\to \mathbb R$ is called to be  countably sub-additive if it satisfies
\begin{description}
  \item[\rm (1)] {\em Countable sub-additivity}: $\Sbep[X]\le \sum_{n=1}^{\infty} \Sbep [X_n]$, whenever $X\le \sum_{n=1}^{\infty}X_n$,
  $X, X_n\in \mathscr{H}$ and
  $X\ge 0, X_n\ge 0$, $n=1,2,\ldots$;
 \end{description}
It is called to be continuous if it satisfies
\begin{description}
  \item[\rm (2) ]  {\em Continuity from below}: $\Sbep[X_n]\uparrow \Sbep[X]$ if $0\le X_n\uparrow X$, where $X_n, X\in \mathscr{H}$;
  \item[\rm (3) ] {\em Continuity from above}: $\Sbep[X_n]\downarrow \Sbep[X]$ if $0\le X_n\downarrow X$, where $X_n, X\in \mathscr{H}$.
\end{description}

\item{\rm (II)}  A function $V:\mathcal F\to [0,1]$ is called to be  countably sub-additive if
$$ V\Big(\bigcup_{n=1}^{\infty} A_n\Big)\le \sum_{n=1}^{\infty}V(A_n) \;\; \forall A_n\in \mathcal F. $$

\item{\rm (III)}  A capacity $V:\mathcal F\to [0,1]$ is called a continuous capacity if it satisfies
\begin{description}
  \item[\rm (III1) ] {\em Continuity from below}: $V(A_n)\uparrow V(A)$ if $A_n\uparrow A$, where $A_n, A\in \mathcal F$;
  \item[\rm (III2) ] {\em Continuity from above}: $V(A_n)\downarrow  V(A)$ if $A_n\downarrow A$, where $A_n, A\in \mathcal F$.
\end{description}
\end{description}
\end{definition}
\bigskip

It is obvious that a continuous sub-additive capacity $V$ (resp. a sub-linear expectation $\Sbep$) is countably sub-additive.
The ``the convergence part'' of the   Borel-Cantelli Lemma is still true for a countably sub-additive capacity.
\begin{lemma} ({\em Borel-Cantelli's Lemma}) Let $\{A_n, n\ge 1\}$ be a sequence of events in $\mathcal F$.
Suppose that $V$ is a countably sub-additive capacity.   If $\sum_{n=1}^{\infty}V\left (A_n\right)<\infty$, then
$$ V\left (A_n\;\; i.o.\right)=0, \;\; \text{ where } \{A_n\;\; i.o.\}=\bigcap_{n=1}^{\infty}\bigcup_{i=n}^{\infty}A_i. $$
\end{lemma}
If $\Capc$ is a continuous capacity and  $\{A_n^c, n\ge 1\}$   are independent relative to $\cCapc$,  i.e., $\cCapc(\bigcap_{j=m}^{m+n} A_j^c)=\prod_{j=m}^n \cCapc(A_j^c)$ for all $n,m\ge 1$, then we can show that $\sum_{n=1}^{\infty}\Capc\left (A_n \right)=\infty$ implies
$ \Capc\left (A_n\;\; i.o.\right)=1$. However, the extended independence does not imply that $\{X_n\in B_n;n\ge 1\}$ are independent relative to $\cCapc$  even when (\ref{eqWeakInd}) is assumed to hold for all non-negative Borel functions $\psi_i$s. So, in general,   we have not ``the divergence part'' of the Borel-Cantelli Lemma.

 Since $\Capc$ may be not countably sub-additive in general, we define an  outer capacity $\outCapc$ by
$$ \outCapc(A)=\inf\Big\{\sum_{n=1}^{\infty}\Capc(A_n): A\subset \bigcup_{n=1}^{\infty}A_n\Big\},\;\; \outcCapc(A)=1-\outCapc(A^c),\;\;\; A\in\mathcal F.$$
  Then it can be shown that $\outCapc(A)$ is a countably sub-additive capacity with $\outCapc(A)\le \Capc(A)$ and the following properties:
  \begin{description}
    \item[\rm (a*)] If $\Capc$ is countably sub-additive, then  $\outCapc\equiv\Capc$;
    \item[\rm (b*)] If $I_A\le g$, $g\in \mathscr{H}$, then $\outCapc(A)\le \Sbep[g]$. Further, if $\Sbep$ is countably sub-additive, then
   $$
    \Sbep[f]\le \outCapc(A)\le \Capc(A)\le \Sbep[g], \;\; \forall f\le I_A\le g, f,g\in \mathscr{H};
  $$
        \item[\rm (c*)]  $\outCapc$ is the largest countably sub-additive capacity satisfying
    the property that $\outCapc(A)\le \Sbep[g]$ whenever $I_A\le g\in \mathscr{H}$, i.e.,
    if $V$ is also a countably sub-additive capacity satisfying $V(A)\le \Sbep[g]$ whenever $I_A\le g\in \mathscr{H}$, then $V(A)\le \outCapc(A)$.
    \end{description}

The following are   our main results on the  Kolmogorov type strong laws of large numbers.

 \begin{theorem} \label{thSLLN1}  Let $\{X_n;n\ge 1\}$ be a sequence identically distributed  random variables in $(\Omega,\mathscr{H},\Sbep)$.
  \begin{description}
  \item[\rm (a)]   Suppose    $\lim\limits_{c\to \infty} \Sbep\left[(|X_1|-c)^+\right]=0$ and $C_{\Capc}[|X_1|]<\infty$. If  $X_1, X_2, \ldots$, are  upper extended negatively  dependent, then
      \begin{equation}\label{eqthSLLN1.0}
\outCapc\left( \limsup_{n\to \infty}\frac{S_n}{n}> \Sbep[X_1] \right)=0.
\end{equation}
If  $X_1, X_2, \ldots$, are    extended negatively  dependent, then
  \begin{equation}\label{eqthSLLN1.1}
\outCapc\left(\Big\{\liminf_{n\to \infty}\frac{S_n}{n}< \cSbep[X_1]\Big\}\bigcup \Big\{\limsup_{n\to \infty}\frac{S_n}{n}> \Sbep[X_1]\Big\}\right)=0.
\end{equation}
\item[\rm (b)]   Suppose that    $\Capc$ is countably sub-additive. Then $\outCapc=\Capc$ and so (a) remains true when $\outCapc$ is replaced by $\Capc$.
\item[\rm (c)] Suppose that  $\Capc$ is continuous.
If $X_1,X_2,\ldots$, are extended independent, and
\begin{equation}\label{eqthSLLN1.3}
\Capc\left(  \limsup_{n\to \infty}\frac{|S_n|}{n}=+\infty \right)<1,
\end{equation}
then $C_{\Capc}[|X_1|]<\infty$.
\end{description}
\end{theorem}

The following corollary shows that the limit of $ {S_n}/{n}$ is a set.
\begin{corollary} \label{cor3.2}
   Let $\{X_n;n\ge 1\}$ be a sequence of extended independent and identically distributed  random variables with
     $C_{\Capc}[|X_1|]<\infty$ and $\lim_{c\to \infty} \Sbep\left[(|X_1|-c)^+\right]=0$. If $\Capc$ is continuous, then
 \begin{equation}\label{eqcor3.2.1}
\Capc\left( \liminf_{n\to \infty}\frac{S_n}{n}= \cSbep[X_1] \right)=1 \;  \text{ and }\;
  \Capc\left(\limsup_{n\to \infty}\frac{S_n}{n}= \Sbep[X_1]\right)=1.
\end{equation}
Moreover, if there is a sequence $\{n_k\}$ with $n_k\to \infty$ and $n_{k-1}/n_k\to 0$ such that $S_{n_{k-1}}$ and $S_{n_k}-S_{n_{k-1}}$ are extended independent, then
 \begin{equation}\label{eqcor3.2.2}
\Capc\left( \liminf_{n\to \infty}\frac{S_n}{n}= \cSbep[X_1]   \;  \text{ and }\;
  \limsup_{n\to \infty}\frac{S_n}{n}= \Sbep[X_1]\right)=1
\end{equation}
and
\begin{equation}\label{eqcor3.2.3}
\Capc\left( C\left\{\frac{S_n}{n}\right\}=\left[\cSbep[X_1], \Sbep[X_1]\right]\right)=1,
\end{equation}
where  $C(\{x_n\})$ denotes the cluster set of a sequence of $\{x_n\}$ in $\mathbb R$.
\end{corollary}

(\ref{eqthSLLN1.1})  tells us that the limit points of $\frac{S_n}{n}$ are between the lower expectation $\cSbep[X_1]$ and the upper
 expectation $\Sbep[X_1]$. (\ref{eqcor3.2.2}) tells us that  the lower expectation   and the upper
 expectation are reachable. (\ref{eqcor3.2.3}) tells us that the interval $\big[\cSbep[X_1], \Sbep[X_1]\big]$ is filled with the limit points.  When $\{X_n;n\ge 1\}$ are independence in the sense Peng's definition, the conclusions in  Theorem \ref{thSLLN1} and Corollary (\ref{cor3.2}) were proved by Zhang (2016). Before that, Chen, Wu and Li (2013) and Chen (2016) proved (\ref{eqthSLLN1.1}) under  a stringer moment condition that $\Sbep[|X_1|^{1+\gamma}]<\infty$ for some $\gamma>0$.

 For a sequence of extended negatively dependent and identically distributed  $\{X_n;n\ge 1\}$ on a probability space $(\Omega,\mathcal F, \pr)$,
 Chen, Chen and Ng (2010) showed that $\pr(S_n/n\to \mu)=1$  if and only if $\ep[|X_1|]<\infty$ and $\ep[X_1]=\mu$.  Under the extended negative dependence in a sub-linear space, we have not find a way to show the conclusions in Theorem \ref{thSLLN1} (c), the inverse part of the strong law of large numbers. However, the conclusions are true if we assume that $\{X_n;n\ge 1\}$ are extended negatively dependent under $\cSbep$ (i.e., in the Definition \ref{Def2.4} $\Sbep$ is replaced by $\cSbep$).

 \begin{theorem}\label{corollary2} Let $\{X_n;n\ge 1\}$ be a sequence of identically distributed  random variables in $(\Omega,\mathscr{H},\Sbep)$ which are extended negatively dependent under $\cSbep$. If
   $\Capc$ is continuous, then
\begin{equation}\label{eqcorollary2.1}
\Capc\left(  \limsup_{n\to \infty}\frac{|S_n|}{n}=+\infty\right)<1 \Longrightarrow C_{\Capc}[|X_1|]<\infty.
\end{equation}
 \end{theorem}

When the sub-linear expectation $\Sbep$ reduces to the linear expectation $\ep$, Theorem \ref{thSLLN1} (b) and Theorem \ref{corollary2} improve the result of   Chen, Chen and Ng (2010).
\begin{corollary}\label{corollary3} Let $\{X_n;n\ge 1\}$ be a sequence of identically distributed  random variables on a probability space $(\Omega,\mathcal{F},\pr)$ which are extended negatively dependent in  the sense of  (\ref{eqENDlower}) and (\ref{eqENDupper}).
If $\ep[|X_1|]<\infty$, then
$\pr\big( S_n/n \to \ep[X_1]\big)=1$.

Conversely,  if $\pr\big(\limsup\limits_{n\to\infty} |S_n|/n=\infty\big)<1$, then $\ep[|X_1|]<\infty$. Further, if $\pr\big(S_n/n\to \mu\big)>0$ for some real $\mu$, then $\ep[|X_1|]<\infty$,  $\mu=\ep[X_1]$ and  $\pr\big( S_n/n\to \mu\big)=1$ .
  \end{corollary}

 According to Corollary \ref{corollary3}, the probability  $\pr\big(S_n/n\to \mu\big)$ is either $0$ or $1$.
\bigskip

For proving the theorems, we need  the following lemma which can be found in Zhang (2016).

\begin{lemma} \label{lem3} Suppose that $X\in \mathscr{H}$ and $C_{\Capc}(|X|)<\infty$.

(a) Then
$$\sum_{j=1}^{\infty} \frac{\Sbep[(|X|\wedge j)^2]}{j^2}<\infty.
$$

(b) Furthermore, if $\lim_{c\to \infty} \Sbep\left[|X|\wedge c\right]=\Sbep \left[|X|\right]$, then
$
\Sbep [|X|]\le C_{\Capc}(|X|).
$
\end{lemma}

{\em Proof of Theorems \ref{thSLLN1}.} We first prove (b).  (a)  follows from (b)
because $\outCapc=\Capc$ when $\Capc$ is countably sub-additive.

It is sufficient to show (\ref{eqthSLLN1.0}) under the assumption that $\{X_n;n\ge 1\}$ are upper extended negatively dependent.
 Without loss of generality, we assume $\Sbep[X_1]=0$. Define $f_c(x)$ and $\widehat{f}(x)$ be defined as in (\ref{eqproofthWLLN1.1})
 and
$\overline{X}_j=f_j(X_j)-\Sbep[f_j(X_j)]$, $\overline{S}_j=\sum_{i=1}^j\overline{X}_i$, $j=1,2,\ldots. $
 Then $f_c(\cdot), \widehat{f}_c(\cdot)\in C_{b,Lip}(\mathbb R)$ and are all non-decreasing functions. And so, $\{\overline{X}_j; j\ge 1\}$, $\{f_j^+(X_j);j\ge  1\}$ and $\{(\widehat{f}_j(X_j))^+;j\ge 1\}$ are all sequences of upper extended negatively   dependent random variables.
  Let $\theta>1$, $n_k=[\theta^k]$. For $n_k<n\le n_{k+1}$, we have
\begin{align*}
 \frac{S_n}{n}= &\frac{1}{n}\left\{ \overline{S}_{n_k}+\sum_{j=1}^{n_k}\Sbep[f_j(X_j)]+\sum_{j=1}^n \widehat{f}_j(X_j)+\sum_{j=n_k+1}^{n}f_j(X_j)\right]\\
 \le & \frac{\overline{S}_{n_k}^+}{n_k}+\frac{\sum_{j=1}^{n_k}|\Sbep[f_j(X_1)]|}{n_k}
 +\frac{ \sum_{j=1}^{n_{k+1}}( \widehat{f}_j(X_j))^+}{n_k}\\
& +\frac{\sum_{j=n_k+1}^{n_{k+1}}\big\{f_j^+(X_j) -\Sbep[f_j^+(X_j)]\big\}}{n_k}
 +\frac{(n_{k+1}-n_k)\Sbep |X_1|}{n_k}\\
=:&(I)_k+(II)_k+(III)_k+(IV)_k+(V)_k.
 \end{align*}
 It is obvious that
 $$\lim_{k\to \infty} (V)_k = (\theta-1)\Sbep [|X_1|] \le (\theta-1)C_{\Capc}(|X_1|)$$
 by Lemma~\ref{lem3} (b).

For $(I)_k$, applying (\ref{eqthIneq2.2ad}) yields
\begin{align*} \Capc\left(\overline{S}_{n_k}\ge \epsilon n_k\right) &
\le C\frac{ \sum_{j=1}^{n_k}\Sbep\big[\overline{X}_j^2\big]}{\epsilon^2 n_k^2}
\le C\frac{  \sum_{j=1}^{n_k}\Sbep\big[f_j^2(X_1)\big]}{\epsilon^2 n_k^2} \\
\le & C\frac{  n_k}{\epsilon^2n_k^2}+C\frac{  \sum_{j=1}^{n_k}\Sbep\big[\big(|X_1|\wedge j)^2\big]}{\epsilon^2 n_k^2}.
\end{align*}
It is obvious that $\sum_k\frac{1}{n_k}<\infty$. Also,
\begin{align*}
 \sum_{k=1}^{\infty}\frac{  \sum_{j=1}^{n_k}\Sbep\big[\big(|X_1|\wedge j)^2\big]}{  n_k^2}
\le &  \sum_{j=1}^{\infty} \Sbep\big[\big(|X_1|\wedge j)^2\}\big]\sum_{k: n_k\ge j}\frac{1}{n_k^2}\\
\le & C\sum_{j=1}^{\infty} \Sbep\big[\big(|X_1|\wedge j)^2\big]\frac{1}{j^2}<\infty
\end{align*}
by Lemma \ref{lem3} (a). Hence
$ \sum_{k=1}^{\infty} \outCapc\left((I)_k\ge \epsilon\right)\le \sum_{k=1}^{\infty} \Capc\left((I)_k\ge \epsilon\right)<\infty.$
 By the Borel-Cantelli lemma and the countable sub-additivity of $\outCapc$, it follows that
$$ \outCapc\left(\limsup_{k\to \infty} (I)_k>\epsilon\right)=0, \;\;\forall \epsilon>0. $$
Similarly,
$$ \outCapc\left(\limsup_{k\to \infty} (IV)_k>\epsilon\right)=0, \;\;\forall \epsilon>0. $$

For $(II)_k$, note that
$$|\Sbep[f_j(X_1)]|=|\Sbep[f_j(X_1)]-\Sbep X_1|\le \Sbep[|\widehat{f}_j(X_1)|]= \Sbep[(|X_1|-j)^+]\to 0. $$
It follows that
$$ (II)_k=\frac{n_{k+1}}{n_k}\frac{\sum_{j=1}^{n_{k+1}}|\Sbep[f_j(X_1)]|}{n_{k+1}}\to 0. $$

At last, we consider $(III)_k$. By the Borel-Cantelli Lemma, we will have
$$ \outCapc\big(\limsup_{k\to \infty}(III)_k>0\big)\le \outCapc\big( \{|X_j|> j\} \; i.o.\big)=0 $$
if we have shown that
\begin{equation}\label{eqproofthSLLN1.10}\sum_{j=1}^{\infty} \outCapc\big(  |X_j|> j \big)\le \sum_{j=1}^{\infty} \Capc\big(  |X_j|> j \big)<\infty.
\end{equation}

Let $g_{\epsilon}$ be a non-decreasing function satisfying that its derivatives of each order are bounded, $g_{\epsilon}(x)=1$ if $x\ge 1$, $g_{\epsilon}(x)=0$ if $x\le 1-\epsilon$, and $0\le g_{\epsilon}(x) \le 1$ for all $x$,
where $0<\epsilon<1$.
Then
\begin{equation} \label{eqfunctionG}
g_{\epsilon}(\cdot)\in C_{b,Lip}(\mathbb R) \text{ is non-decreasing} \;  \text{ and }\; I\{x\ge 1\}\le g_{\epsilon}(x)\le I\{x>1-\epsilon\}.
\end{equation}
Hence by (\ref{eq1.3}),
 \begin{align*}
  \sum_{j=1}^{\infty} \Capc\big(  |X_j|> j \big)
 \le & \sum_{j=1}^{\infty}\Sbep\left[g_{1/2}\big(|X_j|/ j\big)\right]= \sum_{j=1}^{\infty}\Sbep\left[g_{1/2}\big(|X_1|/ j\big)\right]  \;\; (\text{since } X_j\overset{d}= X_1 )\nonumber \\
 \le &  \sum_{j=1}^{\infty} \Capc\big(  |X_1|> j/2 \big)\le 1+C_{\Capc}(2|X_1|)<\infty.\nonumber
 \end{align*}
(\ref{eqproofthSLLN1.10}) is proved. So, we conclude that
$ \outCapc\left(\limsup\limits_{n\to \infty} \frac{S_n}{n}>\epsilon\right)=0$, $\forall \epsilon>0,$
by the arbitrariness of $\theta>1$.  Hence
\begin{align*}
 \outCapc\left(\limsup_{n\to \infty} \frac{S_n}{n}>0\right)
=& \outCapc\left(\bigcup_{k=1}^{\infty}\left\{\limsup_{n\to \infty} \frac{S_n}{n}>\frac{1}{k}\right\}\right)\\
\le &\sum_{k=1}^{\infty}\outCapc\left( \limsup_{n\to \infty} \frac{S_n}{n}>\frac{1}{k} \right)=0.
\end{align*}
(\ref{eqthSLLN1.0}) is proved.

Finally, if $\{X_n; n\ge 1\}$ are lower extended negatively dependent, then
$\{-X_n; n\ge 1\}$ are upper extended negatively dependent. So
$$  \outCapc\left(\liminf_{n\to \infty} \frac{S_n}{n}<\cSbep[X_1]\right)
=\outCapc\left(\limsup_{n\to \infty} \frac{\sum_{k=1}^{n}(-X_k-\Sbep[-X_k])}{n}>0\right)=0. $$
The proof of (\ref{eqthSLLN1.1}) is now completed.

Now, we consider  (c), the inverse part of the strong law of large numbers. Because we have not   ``the divergence part'' of the Borel-Cantelli Lemma and no information about the independence under the conjugate expectation $\cSbep$ or the conjugate capacity $\cCapc$, the proof becomes complex and needs a new approach.  Suppose that $X_1,X_2,\ldots$ are extended independent and identically distributed with $C_{\Capc}(X_1^+)=\infty$. Then, by (\ref{eq1.3}),
\begin{align*}\label{eqproofthSLLN1.11}
 \sum_{j=1}^{\infty}\Sbep\left[g_{1/2}\big(\frac{X_j^+}{Mj}\big)\right]
=&\sum_{j=1}^{\infty}\Sbep\left[g_{1/2}\big(\frac{X_1^+}{Mj}\big)\right]\;\; (\text{ since } X_j\overset{d}= X_1 )\\
\ge &\sum_{j=1}^{\infty} \Capc\big(X_1^+> M j)=\infty, \;\; \forall M>0.\nonumber
\end{align*}
Let $\xi_j=g_{1/2}\big(\frac{X_j^+}{Mj}\big)$, $\eta_n=\sum_{j=1}^n \xi_j$ and $a_n=\sum_{j=1}^n\Sbep[\xi_j]$. Then $a_n\to \infty$ and $\{\xi_j;j\ge 1\}$ are extended independent. For any $0<\delta<\epsilon<1$ and $t>0$, we have
\begin{align*}
&I\left\{\frac{\eta_n-a_n}{a_n}\ge -\epsilon\right\}\ge   e^{-t\delta}\left(\exp\left\{t\frac{\eta_n-a_n}{a_n}\right\}-e^{-t\epsilon}\right)I\left\{\frac{\eta_n-a_n}{a_n}\le  \delta\right\}\\
\ge & e^{-t\delta}\left(\exp\left\{t\frac{\eta_n-a_n}{a_n}\right\}-e^{-t\epsilon}\right)-
e^{-t\delta-t}\exp\left\{t\frac{\eta_n}{a_n}\right\}I\left\{\frac{\eta_n-a_n}{a_n}>  \delta\right\}.
\end{align*}
So
\begin{align*}
 \Capc\left(\frac{\eta_n-a_n}{a_n}\ge -\epsilon\right)
\ge & e^{-t\delta}\left(\Sbep\left[\exp\left\{t\frac{\eta_n-a_n}{a_n}\right\}\right]-e^{-t\epsilon}\right)\\
& -
e^{-t\delta-t}\Sbep\left[\exp\left\{t\frac{\eta_n}{a_n}\right\}g_{1/2}\left(\frac{\eta_n-a_n}{\delta a_n} \right)\right].
\end{align*}
By the extended independence and the fact that $e^x\ge 1+x$, we have
\begin{align*} \Sbep\left[\exp\left\{t\frac{\eta_n-a_n}{a_n}\right\}\right]
=& \prod_{j=1}^n\Sbep\left[\exp\left\{t\frac{\xi_j-\Sbep[\xi_j]}{a_n}\right\}\right]\\
\ge& \prod_{j=1}^n\Sbep\left[ t\frac{\xi_j-\Sbep[\xi_j]}{a_n}+1\right]=1.
\end{align*}
On the other hand, by noting $e^x\le 1+|x|e^{|x|}$, $1+x\le e^x$ and $0\le \xi_j\le 1$,
\begin{align*} \Sbep\left[\exp\left\{2t\frac{\eta_n}{a_n}\right\}\right]
=& \prod_{j=1}^n\Sbep\left[\exp\left\{2t\frac{\xi_j}{a_n}\right\}\right]
\le  \prod_{j=1}^n\Sbep\left[ 1+t\frac{2\xi_j}{a_n}e^{2t/a_n}\right]
\\
\le &\prod_{j=1}^n \left[1+2t\frac{\Sbep[\xi_j]}{a_n}e^{2t/a_n}\right]\le \exp\left\{2t e^{2t/a_n}\right\}.
\end{align*}
Also, by (\ref{eq1.3}),
\begin{align*}
   \Capc\left(\frac{\eta_n-a_n}{a_n}>  \frac{\delta}{2} \right)\le C\frac{4\sum_{j=1}^n \Sbep[(\xi_j-\Sbep[\xi_j])^2]}{\delta^2 a_n^2}
\le C\frac{16\sum_{j=1}^n \Sbep[ \xi_j]}{\delta^2 a_n^2}\le \frac{C}{\delta^2 a_n}.
\end{align*}
It follows that
\begin{align*}
&\Sbep\left[\exp\left\{t\frac{\eta_n}{a_n}\right\}g_{1/2}\left(\frac{\eta_n-a_n}{a_n}\right)\right]
\le   \left\{\Sbep\left[\exp\left\{2t\frac{\eta_n}{a_n}\right\}\right]
\cdot\Sbep\left[g_{1/2}\left(\frac{\eta_n-a_n}{\delta a_n}\right) \right]^2\right\}^{1/2}\\
& \quad \le   \exp\left\{t e^{2t/a_n}\right\}\left\{ \Capc\left(\frac{\eta_n-a_n}{a_n}>  \frac{\delta}{2} \right) \right\}^{1/2} \le  \exp\left\{t e^{2t/a_n}\right\}\frac{C}{\delta a_n^{1/2}} \to 0,
\end{align*}
by  H\"older's inequality and noting $I\{x\ge 1\}\le g_{1/2}(x)\le I\{x\ge 1/2\}$.
We conclude that
$$ \liminf_{n\to \infty}\Capc\left(\frac{\eta_n-a_n}{a_n}\ge -\epsilon\right)
\ge   e^{-t\delta}\left(1-e^{-t\epsilon}\right). $$
Letting $\delta\to 0$ and then $t\to \infty$ yields
\begin{equation}\label{eqproofthSLLN1.12} \lim_{n\to \infty}\Capc\left(\frac{\eta_n-a_n}{a_n}\ge -\epsilon\right)=1.
\end{equation}

Now, choose $\epsilon=1/2$.  By the continuity of $\Capc$,
\begin{align*}
&
\Capc\left(\limsup_{n\to\infty} \frac{X_n^+}{n}>\frac{M}{2}\right)=\Capc\left(\big\{\frac{X_j^+}{Mj}> \frac{1}{2}\big\} \;\;i.o.\right)
\ge
\Capc\left(\sum_{j=1}^{\infty} g_{1/2}\big(\frac{X_j^+}{Mj}\big)=\infty\right)\\
& \qquad = \Capc\left(\Big\{\frac{\eta_n-a_n}{a_n}\ge -\frac{1}{2}\Big\}\;\;i.o. \right)
\ge \limsup_{n\to \infty}  \Capc\left(\frac{\eta_n-a_n}{a_n}\ge -\frac{1}{2}\right)=1.
\end{align*}
On the other hand,
$$\limsup_{n\to\infty} \frac{X_n^+}{n}\le \limsup_{n\to\infty} \frac{|X_n|}{n}\le \limsup_{n\to\infty}\Big(\frac{|S_n|}{n}+\frac{|S_{n-1}|}{n}\Big)\le 2\limsup_{n\to\infty} \frac{|S_n|}{n}.
$$
It follows that
$$ \Capc\left(\limsup_{n\to\infty} \frac{|S_n|}{n}>m\right)=1,\;\; \forall m>0. $$
Hence
$$ \Capc\left(\limsup_{n\to\infty} \frac{|S_n|}{n}=+\infty\right)=\lim_{m\to \infty} \Capc\left(\limsup_{n\to\infty} \frac{|S_n|}{n}>m\right)=1, $$
which  contradicts with (\ref{eqthSLLN1.3}). So, $C_{\Capc}(X_1^+)<\infty$.
Similarly, $C_{\Capc}(X_1^-)<\infty$. It follows that $C_{\Capc}(|X_1|)\le C_{\Capc}(X_1^+)+C_{\Capc}(X_1^-) <\infty$.  \endproof

\bigskip
{\em Proof of Corollary~\ref{cor3.2}.} By (\ref{eqthWLLN1.2}) and the continuity of $\Capc$,
\begin{align*}
  \Capc\left(   \limsup_{n\to \infty}\frac{S_n}{n}\ge  \Sbep[X_1]-\epsilon\right)
\ge   \limsup_{n\to \infty} \Capc\left( \frac{S_n}{n}>\Sbep[X_1]-\epsilon\right)=1, \;\; \forall \epsilon>0.
\end{align*}
B the continuity of $\Capc$ again,
$\Capc\left(   \limsup\limits_{n\to \infty} {S_n}/{n}\ge  \Sbep[X_1] \right)=1, $
which, together with Theorem \ref{thSLLN1} (b)  implies the second equation in (\ref{eqcor3.2.1}). By considering $\{-X_n;n\ge 1\}$ instead, we obtain the first equation in (\ref{eqcor3.2.1}).

For (\ref{eqcor3.2.2}), by noting the facts that  $n_k\to \infty$, $n_{k-1}/n_k\to 0$ such that $S_{n_{k-1}}$ and $S_{n_k}-S_{n_{k-1}}$ are extended independent, we   conclude that
\begin{align*}
&\liminf_{k\to \infty}\Capc\left(\frac{S_{n_{k-1}}}{n_{k-1}}<\cSbep[X_1]+\epsilon \; \text { and }\;  \frac{S_{n_k}-S_{n_{k-1}}}{n_k-n_{k-1}}>\Sbep[X_1]-\epsilon\right) \\
\ge & \liminf_{k\to \infty}\Sbep\left[\phi\left(\frac{S_{n_{k-1}}}{n_{k-1}}-\cSbep[X_1]\right)
\phi\left(\Sbep[X_1]-\frac{S_{n_k}-S_{n_{k-1}}}{n_k-n_{k-1}}\right)\right]\\
\ge & \liminf_{k\to \infty}\Sbep\left[\phi\left(\frac{S_{n_{k-1}}}{n_{k-1}}-\cSbep[X_1]\right)\right]
\cdot \Sbep\left[\phi\left(\Sbep[X_1]-\frac{S_{n_k}-S_{n_{k-1}}}{n_k-n_{k-1}}\right)\right]\\
\ge & \liminf_{k\to \infty}\Capc\left(\frac{S_{n_{k-1}}}{n_{k-1}}<\cSbep[X_1]+\frac{\epsilon}{2}\right)\cdot \Capc\left(  \frac{S_{n_k}-S_{n_{k-1}}}{n_k-n_{k-1}}>\Sbep[X_1]-\frac{\epsilon}{2}\right)\\
 \ge & \liminf_{k\to \infty}\Capc\left(\frac{S_{n_{k-1}}}{n_{k-1}}<\cSbep[X_1]+\frac{\epsilon}{2}\right)\cdot \Capc\left(  \frac{S_{n_k}}{n_k}>\Sbep[X_1]-\frac{\epsilon}{4}\right)= 1, \;\; \forall \epsilon>0,
\end{align*}
by (\ref{eqthWLLN1.1})-(\ref{eqthWLLN1.3}). Hence, by Theorem \ref{thSLLN1} (b) and the continuity of $\Capc$ we have
\begin{align*}
& \Capc\left( \liminf_{n\to \infty}\frac{S_n}{n}\le \cSbep[X_1]+\epsilon  \;  \text{ and }\;
  \limsup_{n\to \infty}\frac{S_n}{n}\ge  \Sbep[X_1]-\epsilon\right) \\
  \ge & \Capc\left( \liminf_{k\to \infty}\frac{S_{n_{k-1}}}{n_{k-1}}\le \cSbep[X_1]+\epsilon  \;  \text{ and }\;
  \limsup_{k\to \infty}\frac{S_{n_k}}{n_k}\ge  \Sbep[X_1]-\epsilon\right) \\
= & \Capc\left( \liminf_{k\to \infty}\frac{S_{n_{k-1}}}{n_{k-1}}< \cSbep[X_1]+\epsilon  \;  \text{ and }\;
  \limsup_{k\to \infty}\frac{S_{n_k}-S_{n_{k-1}}}{n_k-n_{k-1}}>  \Sbep[X_1]-\epsilon\right) \\
  \ge & \Capc\left(\left\{\frac{S_{n_{k-1}}}{n_{k-1}}<\cSbep[X_1]+\epsilon \; \text { and }\;  \frac{S_{n_k}-S_{n_{k-1}}}{n_k-n_{k-1}}>\Sbep[X_1]-\epsilon\right\}\;\; i.o.\right)\\
\ge & \limsup_{k\to \infty} \Capc\left(\frac{S_{n_{k-1}}}{n_{k-1}}<\cSbep[X_1]+\epsilon \; \text { and }\;  \frac{S_{n_k}-S_{n_{k-1}}}{n_k-n_{k-1}}>\Sbep[X_1]-\epsilon\right)=1, \;\; \forall \epsilon>0.
\end{align*}
By the continuity of $\Capc$ again,
$$\Capc\left( \liminf_{n\to \infty}\frac{S_n}{n}\le \cSbep[X_1]   \;  \text{ and }\;
  \limsup_{n\to \infty}\frac{S_n}{n}\ge  \Sbep[X_1] \right)=1, $$
which, together with Theorem \ref{thSLLN1} (b)  implies (\ref{eqcor3.2.2}).

Finally, note
$$ \frac{S_n}{n}-\frac{S_{n-1}}{n-1}=\frac{X_n}{n}-\frac{S_{n-1}}{n-1}\frac{1}{n}\to 0\;\; a.s. \Capc. $$
It can be verified that (\ref{eqcor3.2.2}) implies (\ref{eqcor3.2.3}). \endproof

\bigskip
For proving Theorem \ref{corollary2}, we need  the estimates of $\cCapc\left(S_n\ge x\right)$.

\begin{lemma}\label{lemma2}    Let $\{X_1,\ldots, X_n\}$ be a sequence  of   random variables
in $(\Omega, \mathscr{H}, \Sbep)$ with $\cSbep[X_k]\le 0$ which are upper extended negatively dependent under $\cSbep$ with a dominating constant $K\ge 1$. Then
\begin{description}
  \item[\rm (a)]
  For all $x,y>0$,
$$
\cCapc\left(S_n\ge x\right)\le \Capc\left(\max_{k\le n}X_k\ge y\right)
+ K\exp\left\{-\frac{x^2}{2(xy+B_n)}\Big(1+\frac{2}{3}\ln \big(1+\frac{xy}{B_n}\big)\Big)\right\};
$$
\item[\rm (b)] For any $p\ge 2$, there exists a constant $C_p\ge 1$ such that for all $x>0$ and $0\le \delta\le 1$,
$$
\cCapc\left(S_n\ge x\right)\le C_p K\delta^{-2p} \frac{M_{n,p}}{x^p}
+K\exp\left\{-\frac{x^2}{2 B_n (1+\delta)}\right\};
$$
\item[\rm (c)] We have
\begin{equation}\label{eqlemma2.4}
\cCapc\left(S_n\ge x\right)\le (1+Ke)\frac{\sum_{k=1}^n \Sbep[X_k^2]}{x^2}, \; \forall x>0.
\end{equation}
\end{description}
 \end{lemma}

\proof  Let $Y_k=X_k\wedge y$, $T_n=\sum_{k=1}^n Y_k$ be as  in the proof of Theorem \ref{thIneq2}.   Then
 $$
 \cCapc\left(S_n\ge x\right)\le  \Capc\big(\max_{k\le n}X_k\ge y\big)+
 \cCapc\left(T_n\ge x\right),
 $$
$$
\cCapc\left(T_n\ge x\right)  \le
 e^{-tx}\cSbep[e^{t T_n}]
 \le
 e^{-tx}K\prod_{k=1}^n \cSbep[e^{t Y_k}]
$$
 and
\begin{align*}
\cSbep[e^{tY_k}]\le & \cSbep\left[  1 +tY_k+\frac{e^{ty}-1-t y}{y^2}Y_k^2\right]
\le  1+t\cSbep[Y_k]+\frac{e^{ty}-1-t y}{y^2} \Sbep[Y_k^2]\\
\le & 1+\frac{e^{ty}-1-t y}{y^2} \Sbep[Y_k^2].
\end{align*}
 The remainder proof is similar as that of Theorem \ref{thIneq2}. \endproof

 \bigskip
 {\em The proof of Theorem \ref{corollary2}.} Suppose $C_{\Capc}(X_1^+)=\infty$. Let $g_{\epsilon}(\cdot)$ satisfy (\ref{eqfunctionG}).
 Let $\xi_j=g_{1/2}\big(\frac{X_j^+}{Mj}\big)$, $\eta_n=\sum_{j=1}^n \xi_j$ and $a_n=\sum_{j=1}^n\Sbep[\xi_j]$ be as in the proof of Theorem \ref{thSLLN1} (c). Then $a_n\to \infty$ and $\{-\xi_j-\cSbep[-\xi_j])^2]; j\ge 1\}$ are upper extended negatively dependent under $\cSbep$.
By Lemma \ref{lemma2} (c),
\begin{align*}
&\cCapc\left(\frac{\eta_n-a_n}{a_n}<-\epsilon\right)=  \cCapc\left(\frac{\sum_{j=1}^n (-\xi_j-\cSbep[-\xi_j])}{a_n}>\epsilon\right) \\
\le & (1+Ke) \frac{ \sum_{j=1}^n \Sbep[(-\xi_j-\cSbep[-\xi_j])^2]}{\epsilon^2 a_n^2}\le c \frac{1}{\epsilon^2 a_n}\to 0.
\end{align*}
That is (\ref{eqproofthSLLN1.12}).
By the same argument as in proof of Theorem \ref{thSLLN1}, (\ref{eqproofthSLLN1.12}) will imply a  contradiction to (\ref{eqthSLLN1.3}). So, $C_{\Capc}(X_1^+)<\infty$.
Similarly, $C_{\Capc}(X_1^-)<\infty$. It follows that $C_{\Capc}(|X_1|)<\infty$.  \endproof
  \section{The law of the iterated logarithm}\label{sectLIL}
\setcounter{equation}{0}

In this section, we let $\{X_n;n\ge 1\}$ be a sequence of identically distributed  random variables in $(\Omega,\mathscr{H},\Sbep)$.
 Denote $\overline{\sigma}_1^2=\Sbep[(X_1-\cSbep[X_1])^2]$,
 $\overline{\sigma}_2^2=\Sbep[(X_1-\Sbep[X_1])^2]$,  $a_n=\sqrt{2n\log\log n}$, where $\log x=\ln(x\vee e)$.
The following is the law of the iterated logarithm for   extended negatively dependent random variables.

   \begin{theorem} \label{thLIL}
   Suppose
   $\Sbep[|X_1|^{2+\gamma}]<\infty$ for some $\gamma>0$.
    If $X_1, X_2, \ldots$, are upper extended negatively dependent, then
 \begin{equation}\label{eqthLIL.1}
\outCapc\left(
 \limsup_{n\to \infty}\frac{S_n-n\Sbep[X_1]}{a_n}> \overline{\sigma}_2\right)=0.
\end{equation}
If $X_1, X_2, \ldots$, are extended negatively dependent, then
\begin{equation}\label{eqthLIL.2}
\outCapc\left(\Big\{\liminf_{n\to \infty}\frac{S_n-n\cSbep[X_1]}{a_n}< -\overline{\sigma}_1\Big\}
\bigcup \Big\{\limsup_{n\to \infty}\frac{S_n-n\Sbep[X_1]}{a_n}> \overline{\sigma}_2\Big\}\right)=0.
\end{equation}
\end{theorem}

When the sub-linear expectation $\Sbep$ reduces to the linear expectation, we obtain the law of the iterated logarithm for extended negatively dependent random variables  on a probability space $(\Omega, \mathcal{F}, \pr)$.

\begin{corollary}
   Suppose that $X_1, X_2, \ldots$ are extended negatively dependent and identically distributed  random variables on a probability space $(\Omega, \mathcal{F}, \pr)$  with
   $\ep[|X_1|^{2+\gamma}<\infty$ for some $\gamma>0$ and $\sigma^2=Var(X_1)$.
       Then
\begin{equation}\
\pr\left(\limsup_{n\to \infty}\frac{|S_n-n\ep[X_1]|}{a_n}\le  \sigma \right)=1.
\end{equation}
\end{corollary}

To prove the law of the iterated logarithm, besides the exponential inequality we need a moment inequality on the maximum partial sums.

\begin{lemma}\label{lemMaxIneq}    Let $\{X_k; k=1,\ldots, n\}$ be a sequence  of upper  extended  negatively dependent random variables
in $(\Omega, \mathscr{H}, \Sbep)$ with $\Sbep[X_k]\le 0$, $k=1, \ldots, n$. Let $S_m=\sum_{k=1}^mX_k$, $S_0=0$, $p>2$ (be an integer). And assume that $\Sbep[(|X_k|^p-c)^+]\to 0$ as $c\to \infty$, $k=1, \ldots, n$. Then
\begin{align}\label{eqthmaxIneq.1}
\Sbep\left[\max_{m\le n} (S_m^+)^p\right]\le
C_p n (\log_2 n)^p\max_{k\le n}  C_{\Capc}\Big[(X_k^+)^p\Big]+ C_p n^{p/2} \big(\max_{k\le n} \Sbep[X_k^2]\big)^{p/2}.
\end{align}
 \end{lemma}

\proof We expand $\{X_k; k\le n\}$ to $\{X_k;k\ge 1\}$ by defining $X_k=0$, $k=n+1, n+2,\ldots$.
  Let $T_{k,m}=(X_{k+1}+\cdots+X_{k+m})^+$ and $M_{k,m}=\max_{j\le m}T_{k,j}$. It is easily seen that $T_{k,l+m}\le T_{k,l}+T_{k+l,m}$. Under the conditions in the Lemma, we have $\Sbep[T_{k,m}^p]\le C_{\Capc}(T_{k,m}^p)$. From (\ref{eqthIneq2.3}) it follows that
 $$\Sbep[T_{k,m}^p]\le
C_p m \max_k  C_{\Capc}\Big[(X_k^+)^p\Big]+ C_p m^{p/2} \big(\max_k \Sbep[X_k^2]\big)^{p/2}.
$$
Let $K_1=\Big(C_p\max_k  C_{\Capc}\big[(X_k^+)^p\big]\Big)^{1/p}$ and $K_2= C_p^{1/p}\big(\max_k \Sbep[X_k^2]\big)^{1/2}$. Then
\begin{equation}\label{eqlemmax.1} \Sbep[T_{k,m}^p]\le m \left(K_1+K_2 m^{\frac{p-2}{2p}}\right)^p.
\end{equation}
Using the same argument of M\'orcz (1982), we can show that for some constant $M>1$,
\begin{equation}\label{eqlemmax.2}
\Sbep[M_{k,m}^p]\le M m \left(K_1\log_2 m +K_2 m^{\frac{p-2}{2p}}\right)^p,
\end{equation}
which implies (\ref{eqthmaxIneq.1}).
Here we only give the proof for integer $p$ because it is sufficient for our use. Also, it is sufficient to show that (\ref{eqlemmax.2}) holds for any $m=2^I$.  Suppose (\ref{eqlemmax.2}) holds for $m=2^I$. We will show that it is also true for $m=2^{I+1}$ by the induction.
Now, if $i\le 2^I$, then $T_{k,i}\le M_{k, 2^I}$.  If  $2^I+1\le i\le 2^{I+1}$, then
$T_{k,i}\le T_{k, 2^I}+M_{k+2^I, 2^I}$, and so
$$
T_{k,i}^p\le T_{k, 2^I}^p+M_{k+2^I, 2^I}^p+\sum_{j=1}^{p-1}\binom{p}{j} T_{k, 2^I}^j M_{k+2^I, 2^I}^{p-j}.
$$
It follows that
$$
M_{k,2^{I+1}}^p\le M_{k, 2^I}^p+M_{k+2^I, 2^I}^p+\sum_{j=1}^{p-1}\binom{p}{j} T_{k, 2^I}^j M_{k+2^I, 2^I}^{p-j}.
$$
By the induction,
$$\Sbep\left[M_{k, 2^I}^p\right]\le M 2^I \big(K_1 I+K_22^{\frac{p-2}{2p}I}\big)^{p}, \;\; \Sbep\left[M_{k+2^I, 2^I}^p\right]\le M 2^I \big(K_1 I+K_22^{\frac{p-2}{2p}I}\big)^{p}$$
and
\begin{align*}
 \Sbep\left[ T_{k, 2^I}^j M_{k+2^I, 2^I}^{p-j}\right]
\le & \left\{\Sbep\left[ T_{k, 2^I}^p\right]\right\}^{j/p}\cdot\left\{\Sbep\left[ M_{k+2^I, 2^I}^p\right]\right\}^{(p-j)/p}\\
\le & 2^I  \big(K_1 +K_22^{\frac{p-2}{2p}I}\big)^j\cdot \big\{M^{1/p}\big(K_1 I+K_22^{\frac{p-2}{2p}I}\big)\big\}^{p-j}.
\end{align*}
Let $M>1$ such that $1+M^{-1/p}\le 2^{\frac{p-2}{2p}}$. It follows that
\begin{align*}
&\Sbep\left[M_{k, 2^{I+1}}^p\right]
\le   M 2^I \big(K_1 I+K_22^{\frac{p-2}{2p}I}\big)^{p} \\
&\qquad \qquad +2^I\sum_{j=0}^{p-1} \binom{p}{j}   \big(K_1 +K_22^{\frac{p-2}{2p}I}\big)^j\cdot \big\{M^{1/p}\big(K_1 I+K_2 2^{\frac{p-2}{2p}I}\big)\big\}^{p-j} \\
\le &M 2^I \big(K_1 I+K_2 2^{\frac{p-2}{2p}I}\big)^{p} +2^I\left\{K_1 +K_2 2^{\frac{p-2}{2p}I}+M^{1/p}\big(K_1 I+K_2 2^{\frac{p-2}{2p}I}\big)\right\}^p\\
= &M 2^I \big(K_1 I+K_2 2^{\frac{p-2}{2p}I}\big)^{p} +M 2^I\left\{K_1(I+M^{-1/p}) +K_2 2^{\frac{p-2}{2p}I}(1+M^{-1/p}) \big)\right\}^p\\
\le &2 M 2^I \big(K_1 (I+1)+K_2 2^{\frac{p-2}{2p}(I+1)}\big)^{p}=  M 2^{I+1} \big(K_1 (I+1)+K_2 2^{\frac{p-2}{2p}(I+1)}\big)^{p}.
\end{align*}
The proof is completed. \endproof

\bigskip Now we prove the law of the iterated logarithm.

{\em Proof of Theorem \ref{thLIL}.} It is sufficient to show that (\ref{eqthLIL.1}) under the assumption that $X_1, X_2, \ldots$, are upper extended negatively dependent.
  Without loss of generality, we assume $\Sbep[X_1]=0$ and  $\Sbep[X_1^2]=1$.
 Choose $1/(2+\gamma)<\beta<1/2$, and let $b_n=n^{\beta}$, $a_n=\sqrt{2n\log\log n}$.
Denote $Y_k=(-b_k)\vee (X_k\wedge b_k)$. Then $\{Y_k; k\ge 1\}$ are upper extended negatively dependent.  Note
\begin{align*}
\sum_{n=1}^{\infty}\outCapc(X_n\ne Y_n)\le \sum_{n=1}^{\infty}\Capc(|X_n|>n^{\beta})
\le \sum_{n=1}^{\infty}\frac{\Sbep[|X_1|^{2+\gamma}}{n^{\beta(2+\gamma)}}<\infty.
\end{align*}
Also,
\begin{align*} \sum_{i=1}^n\big| \Sbep[Y_i] \big|=& \sum_{i=1}^n\big| \Sbep[Y_i] -\Sbep[X_i]\big|\le  \sum_{i=1}^n\Sbep[(|X_1|-b_i)^+]\\
\le & \sum_{i=1}^n\frac{\Sbep[|X_1|^{2+\gamma}}{i^{\beta(1+\gamma)}}=O(n^{1-\beta(1+\gamma)})=o(a_n).
\end{align*}
 By the countable  sub-additivity of $\outCapc$,
(\ref{eqthLIL.1}) will follow if we have shown
that
\begin{equation}\label{eqLIL1.12}
\outCapc\left(\limsup_{n\to \infty} \frac{\sum_{k=1}^n (Y_k-\Sbep[Y_k])}{a_n}> (1+\epsilon)^2\right)=0, \;\; \forall \epsilon>0.
\end{equation}

Now, for given $\epsilon$ such that $0<\epsilon<1/2$, let $n_k=[e^{k^{1-\alpha}}]$, where $0<\alpha<\frac{\epsilon}{1+\epsilon}$. Then $n_{k+1}/n_k\to 1$ and $\frac{n_{k+1}-n_k}{n_k}\approx \frac{C}{k^{\alpha}}$. For $n_k<n\le n_{k+1}$,
we have
\begin{align*}
\sum_{i=1}^n (Y_k-\Sbep[Y_k])\le & \sum_{i=1}^{n_k}( Y_i-\Sbep[Y_i])+\max_{n_k<m\le n_{k+1}}\Big(\sum_{i=n_k+1}^m (Y_i-\Sbep[Y_i])\Big)^+
=:  I_k+II_k.
\end{align*}
For the second term, by applying  Lemma \ref{lemMaxIneq}  we have
\begin{align*}
&\Capc\left(II_k\ge  \delta a_{n_k}\right)\le   \frac{\Sbep\left[II_k^p\right]}{(\delta a_{n_k})^p}\\
\le & c \frac{(n_{k+1}-n_k)\big(\log (n_{k+1}-n_k)\big)^p }{a_{n_k}^p} \max\limits_{n_k+1\le i\le n_{k+1}} C_{\Capc}(|Y_i|^p)\\
& +
 c\left( \frac{(n_{k+1}-n_k) }{a_{n_k}^2}\max\limits_{n_k+1\le i\le n_{k+1}}   \Sbep[|Y_i|^2]\right)^{p/2}\\ \le &   c    \frac{(n_{k+1}-n_k)\big(\log (n_{k+1}-n_k)\big)^p }{a_{n_k}^p} n_{k+1}^{\beta(p-2)}  C_{\Capc}(X_1^2)+
 c\left( \frac{n_{k+1}-n_k}{n_k\log\log n_k}\right)^{p/2}\\
 \le &   c    (\log n_k)^p n_k^{-(p-2)(1/2-\beta)}+
 c\left(  \frac{1}{k^{\alpha}}\right)^{p/2}.
 \end{align*}
It follows that
$
\sum_{k=1}^{\infty}\Capc\left(II_k\ge  \delta a_{n_k}\right) <\infty$ for all $\delta>0$
 whenever we choose the integer $p>2$ such that $\alpha p/2>1$. Hence,
\begin{equation}\label{eqLIL1.16} \outCapc\left(\Big\{\frac{II_k}{a_{n_k}}>\delta\Big\}\; i.o.\right)=0,\;\; \forall \delta>0.
\end{equation}

Finally, we consider the   term $I_k$. Let $y=2b_{n_k}$ and $x=(1+\epsilon)^2 a_{n_k}$. Then $|Y_i-\Sbep[Y_i]|\le y$ and
$xy=o(n_k)$. By (\ref{eqthIneq2.1}), we have
\begin{align*}
&\Capc\left(I_k\ge  (1+\epsilon)^2 a_{n_k}\right)\\
\le &
\exp\left\{-\frac{(1+\epsilon)^4 a_{n_k}^2}{2\big(o(n_k)+\sum_{i=1}^{n_k}\Sbep[|Y_i-\Sbep[Y_i]|^2]\big)}\left(1+\frac{2}{3}\ln (1+o(1))\right)\right\}.
\end{align*}
Since
$$ \left|\Sbep[X_i^2]-\Sbep[Y_i^2]\right|\le \Sbep|X_i^2-Y_i^2|=\Sbep[(X_1^2-b_i^2)^+]\to 0, \text{ as } i\to \infty$$
and $|\Sbep[Y_i]|\to 0$  as $i\to \infty$,
we have $\sum_{i=1}^{n_k}\Sbep[|Y_i-\Sbep[Y_i]|^2]\le (1+\epsilon/2)n_k \Sbep X_1^2=(1+\epsilon/2)n_k$ for $k$ large enough.
It follows that
\begin{align*}
  \sum_{k=k_0}^{\infty} & \Capc\left(I_k\ge   (1+\epsilon)^2 a_{n_k}\right)\le
\sum_{k=k_0}^{\infty}\exp\left\{-(1+\epsilon)^2 \log\log n_k\right\}
   \le   \sum_{k=k_0}^{\infty}  \frac{c}{k^{(1+\epsilon)(1-\alpha)}} <\infty
\end{align*}
if $\alpha$ is chosen such that $ (1+\epsilon)(1-\alpha)>1$.  It follows that by the countable sub-additivity and the Borel-Cantelli Lemma again,
\begin{equation}\label{eqLIL1.17} \outCapc\left(\Big\{\frac{I_k}{a_{n_k}}>(1+\epsilon)^2\Big\}\; i.o.\right)=0.
\end{equation}
Combining (\ref{eqLIL1.16}) and (\ref{eqLIL1.17}) yields (\ref{eqLIL1.12}).  The proof is completed. \endproof

\acks
 This work was supported by grants from the NSF of China (No. 11225104), the 973 Program
(No. 2015CB352302) and the Fundamental Research Funds for the Central Universities.



\bigskip

\begin{thebibliography}{99}

{\small \baselineskip 14pt


\bibitem{ChenChenNg10} Chen, Y. Q., Chen, A. Y.  and Ng, K. W. (2010), The strong law of large numbers for extended negatively dependent random variables, {\em Journal of Applied Probability}, {\bf 47}(4): 908-922

\bibitem{Chen16} Chen, Z. J. (2016),  Strong laws of large numbers for sub-linear expectations, {\em Science in China-Mathematics}, {\bf 59}(5): 945-954.  arXiv:1006.0749 [math.PR].

\bibitem{ChenHu14}  Chen, Z. J.  and Hu, F. (2014), A law of the iterated logarithm for sublinear
 expectations, {\em Journal of Financial Engineering}, {\bf 1}, No.02. arXiv: 1103.2965v2[math.PR].

\bibitem{CWL13} Chen, Z. J.,  Wu, P. Y. and  Li, B. M.  (2013), A strong law of large numbers for nonadditive
probabilities, {\em International Journal of Approximate Reasoning}, {\bf 54}:  365-377.

\bibitem{BSS82} Block, H. W., Savits, T. H. and Shaked, M. (1982), Some concepts of negative dependence,
{\em Ann. Probab.}, {\bf 10}: 765-772.


\bibitem{DM06} Denis, L.  and  Martini, C. (2006),  A theoretical framework for the pricing of contingent
claims in the presence of model uncertainty, {\em Ann. Appl. Probab.},  {\bf 16}(2): 827-852.

\bibitem{Gilboa87} Gilboa, I. (1987),  Expected utility theory with purely subjective non-additive prob-
abilities, {\em J. Math. Econom.}, {\bf 16}: 65-68.

\bibitem{JP83} Joag-Dev, K. and Proschan, F. (1983),  Negative association of random variables with
applications, {\em Ann. Statist.}, {\bf  11} (1): 286-295.

\bibitem{Lehmann66} Lehmann, E. (1966), Some concepts of dependence, {\em Ann. Math. Statist.}, {\bf 37} (5):1137-1153.

\bibitem{Liu09} Liu, L. (2009),  Precise large deviations for dependent random variables with heavy tails, {\em Statist. Prob. Lett.}, {\bf 79}: 1290-1298.


\bibitem{MM05} Maccheroni, F. and Marinacci, M. (2005), A strong law of large number for capacities,
{\em Ann. Probab.}, {\bf 33}: 1171-1178.

\bibitem{Mar99}  Marinacci, M. (1999),    Limit laws for non-additive probabilities and their frequentist
interpretation, {\em  J. Econom. Theory}, {\bf 84}: 145-195.

\bibitem{Matula92}  Matula, P.  (1992), A note on the almost sure convergence of sums of negatively dependent
random variables, {\em Statist. Probab. Lett.}, {\bf 15}: 209-213.

\bibitem{Morcz83} M\'oricz (1982), A general moment inequality for the maximum of partial sums of single series, {\em Acta
Sci. Math. Szeged}, {\bf 44}: 67-75.

\bibitem{Newman84} Newman, C. M. (1984), Asymptotic independence and limit theorems for positively and negatively dependent random variables, in {\em Inequalities
in Statistics and Probability} (ed. Tong, Y. L.), IMS Lecture Notes-Monograph Series, Vo1.5,  127-140.


\bibitem{NW81} Newman, C. M.  and Wright, A. L. (1981),  An invariance principle for certain dependent sequences,
{\em Ann. Probab.}, {\bf 9}: 671-675.

\bibitem{Peng99} Peng, S. (1999),  Monotonic limit theorem of BSDE and nonlinear decomposition
theorem of Doob-Meyer type, {\em  Probab. Theory Related Fields},  {\bf 113}: 473-499.

\bibitem{Peng06} Peng, S. (2006),  G-expectation, G-Brownian motion and related stochastic calculus
of Ito type, {\em Proceedings of the 2005 Abel Symposium}.

\bibitem{Peng08a} Peng, S. (2008a), Multi-dimensional G-Brownian motion and related stochastic calculus under G-expectation,
 {\em Stochastic Process. Appl.},  {\bf 118}(12): 2223-2253.

\bibitem{Peng08b} Peng, S. (2008b),  A new central limit theorem under sublinear expectations,
Preprint: arXiv:0803.2656v1 [math.PR]

\bibitem{Peng09} Peng, S. (2009),  Survey on normal distributions, central limit theorem, Brownian
motion and the related stochastic calculus under sublinear expectations, {\em  Science in China Ser. A}, {\bf 52}(7): 1391-1411.

\bibitem{Peng10}  Peng, S. G. (2010), {\em Nonlinear Expectations and Stochastic Calculus under Uncertainty}, arXiv:1002.4546 [math.PR].

\bibitem{Shao00} Shao, Q. M. (2000) A Comparison theorem on moment inequalities between negatively associated and independent random variables,
 {\em J. Theort. Probab.}, {\bf 13}: 343-356.

\bibitem{ShaoSu99}  Shao, Q.M. and Su, C. (1999), The law of the iterated logarithm for negatively associated random variables,
{\em Stochastic Process. Appl.}, {\bf 86}: 139¨C148.

\bibitem{SZW97} Su, C., Zhao, L. C., and Wang, Y. B. (1997), Moment inequalities and weak convergence for negatively associated sequences,
{\em Science in China Ser A}, {\bf 40}: 172-182.

\bibitem{Teran14} Ter\'an, P. (2014),  Laws of large numbers without additivity, {\em Tran.
Amer. Math.  Soc.}, {\bf 366}: 5431-5451.

\bibitem{Zh01a}  Zhang, L. X. (2001a), A Strassen's law of the iterated logarithm for negatively associated random vectors,
{\em  Stoch. Process. Appl.}, {\bf 95}: 311-328

\bibitem{Zh01b} Zhang, L. X. (2001b),
 The weak convergence for functions of negatively associated random variables,  {\em J. Mult. Anal.}, {\bf 78}: 272-298.

 \bibitem{Zh15a} Zhang, L. X. (2015a), Exponential inequalities under sub-linear expectations with applications to laws of the iterated logarithm, Manuscript,  arXiv:1409.0285 [math.PR].

  \bibitem{Zh15b} Zhang, L. X. (2015b), Donsker's invariance principle under the sub-linear expectation with an application to Chung's law of the iterated logarithm, {\em Communications in Math.  Stat.}, {\bf 3}: 187-214. arXiv:1503.02845 [math.PR]


\bibitem{Zh16b} Zhang, L. X. (2016), Rosenthal's inequalities for independent and negatively dependent random variables
  under sub-linear expectations with applications, {\em Science in China-Mathematics},  {\bf 59}(4):751-768.  arXiv:1408.5291 [math.PR].


 }

\end{thebibliography}
\end{document}